%% file: severi-degree.tex
\documentclass{amsart}
\include{headers}

\allowdisplaybreaks

\begin{document}
\title{Severi degrees on toric surfaces}
\author{Fu Liu}
\author{Brian Osserman}
\begin{abstract} 
Ardila and Block used tropical results of Brugall\'e and Mikhalkin to
count nodal curves on a certain family of toric surfaces. Building on a 
linearity result of the first author, we revisit their 
work in the context of the Goettsche-Yau-Zaslow formula for counting nodal
curves on arbitrary smooth surfaces, addressing several 
questions they raised by proving stronger versions of their main theorems. 
In the process, we give new combinatorial formulas for the coefficients
arising in the Goettsche-Yau-Zaslow formulas, and give correction terms
arising from rational double points in the relevant family of toric
surfaces. 
\end{abstract}

\thanks{Fu Liu is partially supported by a grant from the Simons Foundation
\#245939 and by NSF grant DMS-1265702. Brian Osserman is partially supported
by a grant from the Simons Foundation \#279151.}

\maketitle
\tableofcontents

\section{Introduction}

In the 1990's, there was a great deal of work on counting nodal curves in 
surfaces. On the one hand, Ran \cite{ra8}, Kontsevich-Manin
\cite{k-m2} and Caporaso-Harris \cite{c-h1} gave
recursive formulas for counting nodal curves in the plane (generalized
to other rational surfaces by Vakil \cite{va5}), while 
Vainsencher \cite{va6}, Kleiman-Piene \cite{k-p1}, \cite{k-p2} and 
Goettsche \cite{go3} studied the same question on
arbitrary smooth surfaces. The latter line of inquiry, together with
work of Di Francesco and Itzykson \cite{d-i1}, led to remarkable conjectures 
about the structure of formulas counting the curves. We fix notation as 
follows.

\begin{notn}\label{notn:Ndelta}
Let $Y$ be a complex projective surface, and $\sL$ a very
ample line bundle. Given $\delta \geq 0$, let $N^{\delta}(Y,\sL)$ denote
the generalized Severi degree counting the number of curves in $\sL$ having 
$\delta$ nodes, and passing through $\dim |\sL|-\delta$ general fixed 
points of $Y$.
\end{notn}

We will refer to the following as the Goettsche-Yau-Zaslow formula.

\begin{thm}\label{thm:goettsche} For any fixed
$\delta$, there is a polynomial $T_{\delta}(w,x,y,z)$ such that
$$N^{\delta}(Y,\sL)=T_{\delta}(\sL^2,\sL \cdot \sK, \sK^2,c_2)$$
whenever $Y$ is smooth and $\sL$ is sufficiently ample, where 
$\sK$ and $c_2$ are the
canonical class and second Chern class of $Y$, respectively.

Furthermore, if we form the generating function
$$\cN(t)=\sum_{\delta \geq 0} T_{\delta}(w,x,y,z) t^{\delta},$$ 
and set $\cQ(t):= \log \cN(t)$, then 
$$\cQ(t)=w f_1(t)+x f_2(t)+y f_3(t) + z f_4(t)$$
for some $f_1,f_2,f_3,f_4 \in \QQ[[t]]$.

More precisely, there exist power series $B_1(q)$ and $B_2(q)$  such that
$$\sum_{\delta \geq 0}
T_\delta(x, y, z, w)(DG_2(q))^\delta 
= \frac{(DG_2(q)/q)^{\frac{z+w}{12}+\frac{x-y}{2}}
B_1(q)^z B_2(q)^y}
{(\Delta(q)D^2G_2(q)/q^2)^{\frac{z+w}{24}}},$$
where $G_2(q)
= -\frac{1}{24}+ \sum_{n>0}\left(\sum_{d|n} d\right)q^n$
is the second Eisenstein series,
$D=q\frac{d}{dq}$ and $\Delta(q)=q\prod_{k>0}(1-q^k)^{24}$ is
the modular discriminant.
\end{thm}

A symplectic proof of the formula was given by Liu \cite{li4}, and recently 
a simpler algebrogeometric proof was given by Tzeng \cite{tz1}. 
Theorem \ref{thm:goettsche} has 
several distinct components: the first is the universality, asserting that
the answer is given by a polynomial in $\sL^2,\sL \cdot \sK,\sK^2,c_2$.
This was proved also by Kool, Shende and Thomas \cite{k-s-t1},
with a sharper statement on the necessary threshhold of
ampleness for $\sL$; see Remark \ref{rem:ampleness}. The second
component is that if one takes the log of the generating series
$\cN(t)$, one gets linear behavior. We can rephrase this portion of the
theorem as follows.

\begin{notn}\label{notn:Qdelta} Set 
$$Q^{\delta}(Y,\sL)=[t^{\delta}]
\log \left(\sum_{i \geq 0} N^{i}(Y,\sL)t^{i}\right),$$
where $[t^{\delta}] f(t)$ denotes the coefficient of $t^{\delta}$ in
the power series $f(t)$.
\end{notn}

\begin{cor}\label{cor:linearity} For $\sL$ sufficiently ample, we have
$$Q^{\delta}(Y,\sL)=a_{\delta,1} \sL^2 + a_{\delta,2} \sL \cdot \sK 
+ a_{\delta,3} \sK^2 + a_{\delta,4} c_2$$
for some constants $a_{\delta,1},a_{\delta,2},a_{\delta,3},a_{\delta,4}$.
\end{cor}

Finally, of course, there is the specific formula in terms of quasimodular
forms, which generalized a conjecture of Yau and Zaslow for K3 surfaces.
Note that the second and third portions of Theorem \ref{thm:goettsche}
do not refer directly to enumerative geometry, but rather to the structure
of the universal polynomials given in the first part of the theorem.

At the same time as Tzeng's work, Ardila and Block \cite{a-b2} studied the 
same question on a certain family of toric surfaces, which includes many
singular surfaces. The key input for them was the work of Brugalle and
Mikhalkin \cite{b-m3}, which used tropical geometry to combinatorialize the 
problem. A polarized toric surfaces is determined by a polygon
$\Delta$, so we make the following notation.

\begin{notn}\label{notn:NDeltadelta}
Given a lattice polygon $\Delta$, let $(Y(\Delta),\sL)$
be the corresponding polarized toric surface, and set 
$$N^{\Delta,\delta}=N^{\delta}(Y(\Delta),\sL),\quad \text{ and } \quad
Q^{\Delta,\delta}=Q^{\delta}(Y(\Delta),\sL).$$
\end{notn}

A restriction arising from \cite{b-m3} is that we can only
consider polygons which are \textbf{$h$-transverse}, meaning that normal 
vectors have integral (or infinite) slope. Ardila and Block parametrized 
$\Delta$ by vectors $\vec{c},\vec{d}$ describing the normal slopes and 
lattice lengths of the edges, and they showed that if the vertices and
normal directions of $\Delta$ are sufficiently spread out, then
$N^{\Delta,\delta}$ is given by a polynomial in $\vec{c}$ and $\vec{d}$.
This constitutes a universality result, but, as pointed out by the authors,
the need to treat each number of vertices separately, and the lack of 
connection to Theorem \ref{thm:goettsche}, mean that the result is far
from optimal. In the present paper, building on the work \cite{li3}
of the first author, we address three of the four questions 
raised by Ardila and Block in \S 7 of \cite{a-b2}, and simultaneously
give explicit combinatorial formulas for the coefficients $a_{\delta,1},
a_{\delta,2},a_{\delta,3},a_{\delta,4}$ of
Corollary \ref{cor:linearity}. To state our results, we need an additional
definition.

\begin{defn}\label{defn:det}
If $v$ is a vertex of a lattice polygon $\Delta$, we define
the \textbf{determinant}
$\det(v)$ of $v$ to be $\det|w_1,w_2|$, where $w_1$ and $w_2$ are primitive
integer normal vectors to the edges adjacent to $v$.

We say an $h$-transverse polygon $\Delta$ is \textbf{strongly} 
$h$-transverse if either there is a non-zero horizontal edge at the
top of $\Delta$, or the vertex $v$ at the top has $\det(v) \in \{1,2\}$,
and the same holds for the bottom of $\Delta$.
\end{defn}

The condition may appear \textit{ad hoc}, but an $h$-transverse polygon is 
strongly $h$-transverse if and only if the associated toric surface is 
Gorenstein; see Proposition \ref{prop:toric-1}. Corresponding to the 
determinant of a vertex (Proposition \ref{prop:toric-1} again), we define.

\begin{defn} If $Y$ is a surface, a cyclic quotient singularity has
\textbf{index} $d$ if it is analytically isomorphic to $\AA^2/(\ZZ/d\ZZ)$
for some linear action of $\ZZ/d\ZZ$.
\end{defn}

Our main theorem is then the following.

\begin{thm}\label{thm:main} Fix $\delta > 0$. Then there exists an 
(explicitly described) universal linear polynomial
\begin{multline*}
\widehat{T}_\delta(x, y, z, w; s, s_1, \dots,s_{\delta-1} )\\
= a'_{\delta,1}x + a'_{\delta,2}y + a'_{\delta,3}z + a'_{\delta,4}w+
b_{\delta} s + b_{\delta,1}s_1+\dots+b_{\delta,\delta-1} s_{\delta-1}
\end{multline*}
such that if
$\Delta$ is a strongly $h$-transverse polygon with all edges having
length at least $\delta$, and $(Y(\Delta),\sL)$ the corresponding polarized
toric surface, then
$$ Q^{\Delta,\delta} =
\widehat{T}_{\delta}(\sL^2, \sL \cdot \sK,
\sK^2, \tc_2, S, S_1,\dots,S_{\delta-1}),$$
where $\sK$ is the canonical line bundle on $Y(\Delta)$,
$S_i$ is the number of singularities of $Y(\Delta)$ of index $i+1$,
$\tc_2=c_2(Y(\Delta))+\sum_{i \geq 1} i S_i$, 
and $S=\sum_{i \geq 1} (i+1) S_i$.
\end{thm}

Note that $\tc_2=c_2$ in the nonsingular case. In general (under
the strongly $h$-transverse hypothesis),
$\tc_2$ is the second Chern class of a minimal desingularization
of $Y(\Delta)$.

In fact, our calculation treats arbitrary $h$-transverse polygons,
with additional correction terms for non-Gorenstein singularities.
However, we have yet to find an equally satisfactory
formulation in this case;
see Theorem \ref{thm:main-2} for this as well as the precise formulas for
the coefficients of $\widehat{T}_{\delta}$.

The proof of Theorem \ref{thm:main} is independent of Theorem 
\ref{thm:goettsche}, but a simple argument shows the following 
compatibility.

\begin{prop}\label{prop:coefs-agree} The coefficients $a'_{\delta,i}$
arising in Theorem \ref{thm:main} agree with the corresponding 
coefficients $a_{\delta,i}$ of Corollary \ref{cor:linearity}.
\end{prop}

We can then conclude the following.

\begin{cor}\label{cor:goettsche-corrected}
If we set
\[ T_\delta(x, y, z, w; s, s_1, \dots s_{\delta-1} ) 
:= [t^\delta] \exp \left( \sum_{i \ge 1} 
\widehat{T}_i(x, y, z, w; s, s_1, \dots, s_{i-1} ) t^\delta \right),\]
we have 
\begin{multline*}
   \sum_{\delta \geq 0}T_\delta(x, y, z, w; s, s_1, \dots, s_{\delta-1})(DG_2(q))^\delta \\
  =  \frac{(DG_2(q)/q)^{\frac{z+w}{12}+\frac{x-y}{2}}B_1(q)^{z}B_2(q)^{y}}{(\Delta(q)D^2G_2(q)/q^2)^{\frac{z+w}{24}}} \cP(q)^{-s} \prod_{i \ge 2} \cP\left(q^i \right)^{s_{i-1}},
\end{multline*}
where $\cP(x) = \sum_{n \ge 0} p(n) x^n$ is the generating function for
partitions, and other notation is as in Theorem \ref{thm:goettsche}.
\end{cor}

Thus, using Theorem \ref{thm:main}, we obtain natural correction factors
to the formulas of Theorem \ref{thm:goettsche} arising from the presence
of rational double points, and we show that these factors give the correct 
enumerative formulas
in the case of toric surfaces arising from strongly $h$-transverse
polygons. At the same time, we give new explicit combinatorial formulas
for the coefficients appearing in Corollary \ref{cor:linearity}, which
are known to give valid enumerative formulas for all smooth surfaces.
In particular, we give combinatorial (but very complicated) formulas
for the power series $B_1(q),B_2(q)$ appearing in Theorem 
\ref{thm:goettsche}; see Corollary \ref{cor:b1-b2} below.

We briefly discuss our techniques. We start from two basic ingredients:
first, the 
combinatorial formula for $N^{\Delta,\delta}$ deduced by Ardila and
Block from the work of Brugalle and Mikhalkin, which we review in
\S \ref{sec:bg-severi-deg} below, and second, a linearity results which
is the main theorem of \cite{li3}, and discussed in \S 
\ref{sec:bg-linearity}. The aforementioned formula for $N^{\Delta,\delta}$ 
is expressed as a sum over
certain graphs, called ``long-edge graphs.'' The sum is further complicated
by the need to consider reorderings of the ``left and right edge directions'' 
of the polygon $\Delta$. Instead of attempting to 
analyze $N^{\Delta,\delta}$ directly, as Ardila and Block did, following
the recent work of Block, Colley and Kennedy \cite{b-c-k1} we 
instead focus our attention throughout on $Q^{\Delta,\delta}$. The main
theorem of \cite{li3}, due to the first author, is that in this context
the functions of interest become linear, at least within certain ranges. 
This result was originally motivated by a conjecture of Block, Colley and 
Kennedy in the special case of counting plane curves, but the proof goes 
through in the more general setting needed for the present paper.

With the above preliminaries out of the way, in \S
\ref{sec:deviation} we analyze how far the actual formulas deviate from
the linear functions described in \S \ref{sec:bg-linearity}, and in
\S \ref{sec:no-reorder} we derive a formula for the value of the relevant
sums of the linear functions in question, for the default ordering of the
edge directions of $\Delta$. In \S \ref{sec:reorder} we consider how this
formula varies with the choice of reordering, obtaining our first results
in the direction of Theorem \ref{thm:main}. Finally, in \S \ref{sec:geometry}
(supplemented by calculations on toric surfaces from Appendix \ref{app:toric}),
we bring geometry into the picture, using some tricks arising from the
instrinsic nature of $N^{\Delta,\delta}$ to simplify our formulas and
derive our main results. We conclude in \S \ref{sec:examples} with
some discussion of examples for small values of $\delta$.

\subsection*{Acknowledgements} We would like to thank Florian Block,
Erwan Brugall\'e and Steven Kleiman for helpful conversations.

\section{Severi degrees via long-edge graphs}\label{sec:bg-severi-deg}

We begin by reviewing the combinatorial objects which arise in the
formula of Brugalle-Mikhalkin and Ardila-Block for Severi degrees.
The basic objects are the following.

\begin{defn}\label{defn:long-edge}
A \textbf{long-edge graph} $G$ is a graph $(V,E)$ with a weight function 
$\rho$ satisfying the following conditions:
\begin{alist}
\itm The vertex set is fixed as $V = \{0, 1, 2, \dots\},$ and the edge set 
$E$ is finite.
\itm Multiple edges are allowed, but loops are not.
\itm The weight function $\rho: E \to \ZZ_{>0}$ assigns a positive integer 
to each edge.
\itm There is no ``short edge,'' i.e., there is no edge connecting $i$ and 
$i+1$ with weight $1.$
\end{alist}
\end{defn}

If an edge of $G$ connects $i$ to $j$ with $i<j$, we will say it is
``from $i$ to $j$.''

Two fundamental numbers associated to long-edge graphs are as follows.

\begin{defn}\label{defn:mudelta}
Given a long-edge graph $G = (V, E)$ equipped with weight function $\rho,$ 
we define the \textbf{multiplicity of $G$} to be 
\[ \mu(G) = \prod_{e \in E} (\rho(e))^2,\]
and the \textbf{cogenus of $G$} to be
\[ \delta(G) = \sum_{e \in E} \left( \ell(e) \rho(e) - 1 \right),\]
where for any $e \in E$ from $i$ to $j$, we define 
$\ell(e) = j-i.$
Note that any non-empty long-edge graph has positive cogenus.
\end{defn}

\begin{figure}
                \begin{picture}(440,60)(-15,-15)
\setlength{\unitlength}{4.0pt}\thicklines
\multiput(0,2)(8,0){3}{\circle*{1}}
\put(0,2){\line(1,0){8}}
\qbezier(0,2)(8,12)(16,2)
\put(4,0){\makebox(0,0){$\scriptstyle 2$}}
\put(8,9){\makebox(0,0){$\scriptstyle 1$}}
\put(0,-0.5){\makebox(0,0){$0$}}
\put(8,-0.5){\makebox(0,0){$1$}}
\put(16,-0.5){\makebox(0,0){$2$}}
\put(8,-4.5){\makebox(0,0){$G_1$}}

\multiput(30,2)(8,0){3}{\circle*{1}}
\put(30,2){\line(1,0){8}}
\qbezier(30,2)(38,12)(46,2)
\put(34,0){\makebox(0,0){$\scriptstyle 2$}}
\put(38,9){\makebox(0,0){$\scriptstyle 1$}}
\put(30,-0.5){\makebox(0,0){$3$}}
\put(38,-0.5){\makebox(0,0){$4$}}
\put(46,-0.5){\makebox(0,0){$5$}}
\put(38,-4.5){\makebox(0,0){$G_2$}}

\multiput(60,2)(8,0){4}{\circle*{1}}
\put(60,2){\line(1,0){8}}
\qbezier(60,2)(68,12)(76,2)
\put(76,2){\line(1,0){8}}
\put(64,0){\makebox(0,0){$\scriptstyle 2$}}
\put(68,9){\makebox(0,0){$\scriptstyle 1$}}
\put(80,0){\makebox(0,0){$\scriptstyle 2$}}
\put(60,-0.5){\makebox(0,0){$3$}}
\put(68,-0.5){\makebox(0,0){$4$}}
\put(76,-0.5){\makebox(0,0){$5$}}
\put(84,-0.5){\makebox(0,0){$6$}}
\put(72,-4.5){\makebox(0,0){$G_3$}}
\end{picture}
\caption{Examples of long-edge graphs}
\label{fig:exlongedge}
\end{figure}
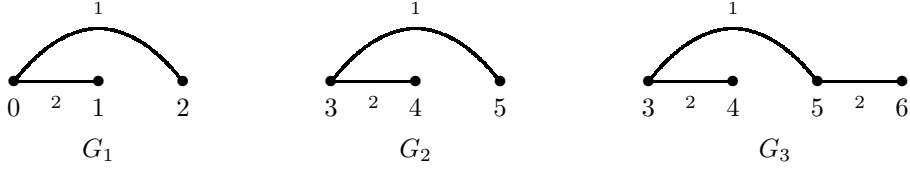

\begin{ex}\label{ex:longedge1}
Figure \ref{fig:exlongedge} shows three long-edge graphs. We have
$$\mu(G_1)=\mu(G_2)=4, \quad \text{ and } \quad \delta(G_1)=\delta(G_2)=2.$$
We also have
$$\mu(G_3)=16, \quad \text{ and } \quad \delta(G_3)=3.$$
\end{ex}

The following definitions will also be important.

\begin{defn}\label{defn:shift}
We define $\minv(G)$ (respectively, $\maxv(G)$) to be the smallest 
(respectively, largest) vertex of $G$ that has nonzero degree. We then 
define the \textbf{length} of $G,$ denoted by $\ell(G),$ to be 
$\maxv(G) -\minv(G).$

For any long-edge graph $G$ and any $k \in \ZZ_{\geq 0},$ we denote by 
$G_{(k)}$ the 
graph obtained by shifting all edges of $G$ to the right $k$ units, i.e., 
a weighted edge from $i$ to $j$ in $G$ becomes a weighted edge from
$i+k$ to $j+k$ in $G_{(k)}.$
\end{defn}

The formula for Severi degrees is expressed as a sum over certain long-edge
graphs of a rather complicated associated statistic. A preliminary 
definition is the following. 

\begin{defn}\label{defn:lambda}
Let $G$ be a long-edge graph with associated weight function $\rho$. For any
$j$, we define
\[ \lambda_j(G) = 
\sum_{e} \rho(e), \]
where $e$ ranges over edges of $G$ 
from $i$ to $k$ with $i<j \leq k$.
Let $\bbeta= (\beta_0, \beta_1, \dots, \beta_{M}) \in \ZZ_{\geq 0}^{M+1}$ 
(where $M \ge 0$).
We say $G$ is \textbf{$\bbeta$-allowable} if $\maxv(G) \le M+1$ and 
$\beta_{j-1} \ge \lambda_j(G)$ for each $j.$

A long-edge graph $G$ is \textbf{strictly $\bbeta$-allowable} if it satisfies 
the following conditions:
\begin{alist}
\itm $G$ is $\bbeta$-allowable.
\itm Any edge that is incident to the vertex $0$ has weight $1.$ 
\itm Any edge that is incident to the vertex $M+1$ has weight $1.$ 
\end{alist}
\end{defn}

The relevant statistic is then the $P^s_{\bbeta}(G)$ defined below.
In fact, its precise definition will not be relevant to us, but we
recall it for the sake of completeness.

\begin{defn}\label{defn:P}
Suppose $\bbeta= (\beta_0, \beta_1, \dots, \beta_{M}) \in \ZZ_{\geq 0}^{M+1}$ 
and $G$ is $\bbeta$-allowable. We create a new graph $\ext_{\bbeta}(G)$ by 
adding $\beta_{j-1}-\lambda_j(G)$ (unweighted) edges connecting vertices 
$j-1$ and $j$ for each $1 \le j \le M+1.$ 

An \textbf{$\bbeta$-extended ordering} of $G$ is a total ordering of the 
vertices and edges of $\ext_{\bbeta}(G)$ satisfying the following conditions:
\begin{alist}
\itm The ordering extends the natural ordering of the vertices 
$0, 1, 2, \cdots$ of $\ext_{\bbeta}(G)$.
\itm For any from $a$ to $b$,
its position in the total ordering has to be between $a$ and $b.$
\end{alist}

We consider two extended orderings $o$ and $o'$ to be equivalent if there is 
an automorphism $\sigma$ on the edges of $\ext_{\bbeta}(G)$ such that
\begin{alist}
\itm If $\sigma(e) = e',$ then $e$ and $e'$ have the same vertices, and 
either have the same weights or are both unweighted.
\itm When applying $\sigma$ on the ordering $o$, one obtains the ordering 
$o'.$
\end{alist}
For any long-edge graph $G$, we define
\[ P_{\bbeta}(G) = \text{ the number of $\bbeta$-extended orderings 
(up to equivalence) of $G$},\]
where by convention $P_{\bbeta}(G) = 0$ if $G$ is not $\bbeta$-allowable, and 
then define
\[ P_{\bbeta}^s(G) = \begin{cases}P_{\bbeta}(G) & \text{if $G$ is strictly 
$\bbeta$-allowable;} \\
	0 & \text{otherwise.}
\end{cases}\]
\end{defn}

\begin{rem}\label{rem:emptygraph}
Suppose $G$ is the empty long-edge graph, i.e., the graph without any edges. 
Then for any $\bbeta \in \ZZ_{\geq 0}^{M+1},$ the graph $G$ is (strictly) 
$\bbeta$-allowable and $P_{\bbeta}(G) = P_{\bbeta}^s(G) =\mu(G)= 1.$
\end{rem}

We will frequently find it convenient to specify $\bbeta$ via its
successive differences, so we 
set the following notation.

\begin{notn}\label{notn:divergence}
Given a 
sequence $\bd = (d_0, d_1, \dots, d_M),$ we set
\[ \bbeta(\bd) = (d_0, d_0+d_1, d_0+d_1 + d_2, \dots, d_0 + d_1 + \cdots 
+ d_M).\]
\end{notn}

In order to describe which long-edge graphs to sum over, we
need the following definitions.

\begin{defn} Given an $h$-transverse polygon $\Delta$, 
The \textbf{left} (respectively,
\textbf{right}) \textbf{directions} of $\Delta$ are defined to be the
(elements of the) integer multiset consisting of the slopes of the 
normal vectors at each half-integer height on the left (respectively, 
right) side of $\Delta$.
\end{defn}

Thus, if $\Delta$ has height $M$, then each of the left and right directions
consist of $M$ elements, if counted with multiplicity.

\begin{defn}\label{defn:Rev}
The \textbf{reversal set} $\Rev(\bs)$ of a sequence $\bs=(s_1, \dots, s_M)$ is 
\[ \Rev(\bs) = \{1 \le i < j \le M: s_i < s_j \}.\]

Therefore, if $\bl$ and $\br$ are reorderings of the multisets of left and
right directions of $\Delta$, the reversal sets of $\br$ and $-\bl$ are
\[
\Rev(\br) = \{1 \le i < j \le M: r_i < r_j \}, \quad \Rev(-\bl) 
= \{1 \le i < j \le M: l_i > l_j \}.
\]
Define the \textbf{cogenus} of the pair $(\bl, \br)$ to be
\[
\delta(\bl, \br) = \sum_{(i,j) \in \Rev(\br)} (r_j - r_i) \quad 
+ \sum_{(i,j) \in \Rev(-\bl)} (l_i - l_j).
\]
\end{defn}

We also use the following notation.

\begin{notn}\label{notn:dt}
Given an $h$-transverse polygon $\Delta$,
let $d^t$ (respectively, $d^b$) denote the width of $\Delta$ at the top 
(respectively, bottom). 
\end{notn}

Finally, we can describe the formula for the Severi 
degree using long-edge graphs, as deduced in \cite{a-b2} from
\cite{b-m3}. 

\begin{thm}[\cite{a-b2}, Proposition 3.3]
For any $\delta \geq 0$, and any $h$-transverse polygon $\Delta$ with
all sides of length at least $\delta-1$,
the Severi degree $N^{\Delta,\delta}$ is given by
$$N^{\Delta,\delta} 
= \sum_{(\bl, \br)} \ \sum_{G} \mu(G) P_{\bbeta(d^t, \br-\bl)}^s (G).$$
where the first summation is over all reorderings $\bl=(l_1, \dots, l_M)$ 
and $\br=(r_1, \dots, r_M)$ of the multisets of left and 
right directions of $\Delta$ satisfying $\delta(\bl,\br) \le \delta$ and
$\bbeta(d^t,\br-\bl) \in \ZZ_{\geq 0}^{M+1}$, and 
the second summation is over all long-edge graphs $G$ of cogenus 
$\delta - \delta(\bl, \br).$
\end{thm}

In principle, it is possible for a choice of $(\bl,\br)$ to yield a
sequence $\bbeta(d^t,\bl-\br)$ for which some entries are negative. 
However, for such choices it is natural to set 
$P_{\bbeta(d^t,\br-\bl)}^s(G)=0$, so we will always require that
$\bbeta(d^t,\br-\bl) \in \ZZ_{\geq 0}^{M+1}$.

Compared to the phrasing in \cite{a-b2}, the need for the sides of
$\Delta$ to have length at least $\delta-1$ is due to the phrasing 
of Severi degrees in terms of counting curves in $Y(\Delta)$ as opposed
to curves in the torus of prescribed Newton polygon; see the discussion
immediately prior to the statement of Theorem 1.2 in \cite{a-b2}.

\section{Templates and linearity results}\label{sec:bg-linearity}

In this section, we shift our attention from $N^{\Delta,\delta}$ to
$Q^{\Delta,\delta}$, and recall the resulting simplifications, most
importantly the linearity result (Theorem \ref{thm:linear0} below) from
\cite{li3}.

Given $\bbeta \in \ZZ_{\geq 0}^{M+1}$, define
\begin{equation}\label{equ:NbetaDelta}
N_{\bbeta}^{\delta} = \sum_G \mu(G) P_{\bbeta}^s(G), 
\end{equation}
where the summation is over all long-edge graphs $G$ of cogenus $\delta.$ Then
\[ N^{\Delta, \delta} = \sum_{(\bl, \br)} N_{\bbeta(d^t, \br-\bl)}
^{\delta - \delta(\bl, \br)},\]
where the summation is over all choices of $(\bl,\br)$ with
$\bbeta(d^t,\br-\bl) \in \ZZ_{\geq 0}^{M+1}$ and $\delta(\bl,\br) \leq \delta$.

Consider the generating functions 
\begin{equation}\label{equ:NDeltat}
\cN^{\Delta}(t) := 1 + \sum_{\delta \ge 1} N^{\Delta, \delta} t^{\delta}
\end{equation}
and
\begin{equation}\label{equ:Nbetat}
\cN(\bbeta,t) := 1 + \sum_{\delta \ge 1} N_{\bbeta}^{\delta} t^{\delta}
= 1 + \sum_{G} \mu(G) P_{\bbeta}^s(G) t^{\delta(G)},
\end{equation}
where the last sum above is over non-empty long-edge graphs $G$.
Then
\[ \cN^{\Delta}(t) = \sum_{(\bl, \br)} \cN(\bbeta(d^t, \br-\bl),t) 
t^{\delta(\bl, \br)},\]
where the summation is over all choices of $(\bl,\br)$ with
$\bbeta(d^t,\br-\bl) \in \ZZ_{\geq 0}^{M+1}$.
Let 
\begin{equation}\label{equ:QDeltat}
\cQ^{\Delta}(t) :=  \sum_{\delta \ge 1} Q^{\Delta, \delta} t^{\delta}
\end{equation}
be the formal logarithm of $\cN^{\Delta}(t)$,
and also let $\cQ(\bbeta,t)$ be the formal logarithm of $\cN(\bbeta,t)$ and 
$Q_{\bbeta}^{\delta}$ be the coefficient of $t^{\delta}$ in $\cQ(\bbeta,t),$ 
i.e.,
\begin{equation}\label{equ:Qbetat}
 \cQ(\bbeta,t) := \log \left(\cN(\bbeta,t)  \right) 
= \sum_{\delta \ge 1} Q_{\bbeta}^{\delta} t^{\delta}.
\end{equation}
Then
\begin{equation}
Q_{\bbeta}^{\delta} = \sum_{i \ge 1} \frac{(-1)^{i+1}}{i} 
\sum_{(G_1, \dots, G_i)}  \left(\prod_{j=1}^i \mu(G_j) P_{\bbeta}^s(G_j)\right),
\label{equ:Qbeta0}
\end{equation}
where the summation is over all the tuples $(G_1, \dots, G_i)$ of 
(non-empty) long-edge graphs satisfying $\sum_{j=1}^i \delta(G_j) = \delta.$

Because the constant coefficients of $\cQ^{\Delta}(t)$ and $\cQ(\bbeta,t)$ 
vanish, and we will mainly be studying these generating series, our
convention throughout will be that $\delta \geq 1$ unless we specify
otherwise.

\begin{defn}
Given a long-edge graph $G,$ we say a tuple $(G_1, \dots, G_i)$ of 
(non-empty) long-edge graphs is a \textbf{partition} of $G$ if the disjoint 
union of the (weighted) edge sets of $G_1, \dots, G_i$ is the (weighted) 
edge set of $G.$
\end{defn}

By the definitions of multiplicity and cogenus, one checks that for any 
partition $(G_1, \dots, G_i)$ of $G$, we have
\begin{equation}\label{equ:partofG}
\mu(G) = \prod_{j=1}^i \mu(G_j) \quad \text{and} \quad 
\delta(G) = \sum_{j=1}^i \delta(G_j).
\end{equation}

Since we can consider any such tuple $(G_1,\dots,G_i)$ a partition of a 
long-edge graph of 
cogenus $\delta,$ it is natural to give the following definition: for any 
long-edge graph $G$, we define
\begin{equation}\label{equ:phibetas}
\Phi_{\bbeta}^s(G) := \sum_{i \ge 1} \frac{(-1)^{i+1}}{i} 
\sum_{(G_1, \dots, G_i)}  \left(\prod_{j=1}^i P_{\bbeta}^s(G_j) \right),
\end{equation}
where the summation is over all the partitions of $G$.

With this definition and by \eqref{equ:partofG}, we can rewrite 
\eqref{equ:Qbeta0}:
\begin{equation}\label{equ:Qbeta} 
	Q_{\bbeta}^{\delta} = \sum_G \mu(G) \Phi_{\bbeta}^s(G),
  \end{equation}
where the summation is over all the long-edge graphs of cogenus $\delta$. 

We define $\Phi_{\bbeta}(G)$ analogously:
\begin{equation}\label{equ:phibeta}
\Phi_{\bbeta}(G) := \sum_{i \ge 1} \frac{(-1)^{i+1}}{i} \sum_{(G_1, \dots, G_i)}  \left(\prod_{j=1}^i P_{\bbeta}(G_j) \right).
\end{equation}

\begin{defn}
A long-edge graph $\Gamma$ is a \textbf{template} if for any vertex 
$i:$ $1 \le i \le \ell(\Gamma)-1$, there exists at least one edge from $j$
to $k$ satisfying $j < i < k.$

We say a long-edge graph $G$ is a \textbf{shifted template}  if $G$ can be obtained by shifting a template; that is, if $G = \Gamma_{(k)}$ for some template $\Gamma$ and some nonnegative integer $k.$
\end{defn}

\begin{ex}\label{ex:longedge2} Continuing with Example \ref{ex:longedge1},
we see that $G_1$ is a template, $G_2$ is a shifted template but not a 
template, and $G_3$ is not a shifted template.
\end{ex}

It is not difficult to see that only shifted templates contribute to the 
``logarithmic'' version Severi degrees. That is, we have the following.

\begin{lem}\label{lem:vanish0}
Suppose $G$ is not a shifted template. Then
\[ \Phi_{\bbeta}^s(G) = 0.\]
\end{lem}

See for instance Lemma 2.15 of \cite{li3}.

\begin{notn}\label{notn:epsilon}
Let $G$ be a long-edge graph. We define
$$ \epsilon_0(G) = \begin{cases}1, 
& \text{if all edges adjacent to the vertex $\minv(G)$ have weight $1$;} \\
0, & \text{otherwise.}\end{cases}$$
$$ \epsilon_1(G) = \begin{cases}1, & \text{if all edges 
adjacent to the vertex $\maxv(G)$ have weight $1$;} \\
0, & \text{otherwise.}\end{cases}
$$
\end{notn}

Putting together \eqref{equ:Qbeta} with Lemma \ref{lem:vanish0} and
Corollary 3.5 of \cite{li3}, we conclude the following.

\begin{cor}\label{cor:Qbetaformula}
  \[ Q_{\bbeta}^{\delta} = \sum_{\Gamma} \mu(\Gamma) \sum_{k=1-\epsilon_0(\Gamma)}^{M-\ell(\Gamma) + \epsilon_1(\Gamma)} \Phi_\bbeta \left( \Gamma_{(k)} \right),\]
where the first summation is over all templates $\Gamma$ of cogenus $\delta.$
\end{cor}

One more definition will be important to us.

\begin{defn}\label{defn:olam}
For a long-edge graph $G$, and any $j \geq 1$, define
\[\olam_j(G):=\lambda_j(G)-\#\{\text{edges from $j-1$ to $j$}\}.\]
We say $G$ is \textbf{$\bbeta$-semiallowable} if $\maxv(G) \le M+1$ and 
$\beta_{j-1} \ge \olam_j(G)$ for each $j.$
\end{defn}

The following is the main result (Theorem 1.4) of \cite{li3}, and 
together with Corollary \ref{cor:Qbetaformula} constitutes the starting
point for our work.

\begin{thm}\label{thm:linear0} Suppose $G$ is a long-edge graph of length 
$\ell.$ There exists a linear multivariate function $\Phi(G, \bbeta)$ in 
$\bbeta$ such that for any 
$\bbeta$ with $G$ is $\bbeta$-semiallowable, we have
$\Phi_{\bbeta}(G)=\Phi(G,\bbeta)$.
\end{thm}

\begin{rem}\label{rem:linear-support}
It is immediate from the definition that $P_{\bbeta}(G)$,
and hence $\Phi_{\bbeta}(G)$, only depends on the value of $\beta_i$ 
for $i$ with $\minv(G) \leq i \leq \maxv(G)-1$. It then follows that the 
linear function $\Phi(G,\bbeta)$ of Theorem \ref{thm:linear0} likewise
has nonzero coefficient for $\beta_i$ only if 
$\minv(G) \leq i \leq \maxv(G)-1$.
\end{rem}

Theorem \ref{thm:linear0} allows us to define the following ``linearized'' 
version of $Q_{\bbeta}^{\delta}$.

\begin{notn}\label{notn:QGammabbeta}
  Given $\bbeta$ and a template $\Gamma,$ define
  \[ Q_\Gamma(\bbeta) = \mu(\Gamma) \sum_{k=1-\epsilon_0(\Gamma)}^{M-\ell(\Gamma) + \epsilon_1(\Gamma)} \Phi \left( \Gamma_{(k)}, \bbeta \right).\]
  Then define
  \[ Q_{\delta}(\bbeta) = \sum_{\Gamma} Q_\Gamma(\bbeta) = \sum_{\Gamma} \mu(\Gamma) \sum_{k=1-\epsilon_0(\Gamma)}^{M-\ell(\Gamma) + \epsilon_1(\Gamma)} \Phi \left( \Gamma_{(k)}, \bbeta \right),\]
where the summations indexed by $\Gamma$ range over all templates of cogenus 
$\delta.$
\end{notn}

We also introduce notation for the coefficients and important related
combinations.

\begin{notn}\label{notn:zetaeta}
Suppose $\Gamma$ is a template of length $\ell,$ and 
$\Phi(\Gamma, \bbeta) = \eta_0 + \sum_{j=0}^{\ell-1} \eta_{j+1} \beta_j.$ We define
\[ \zeta^{i}(\Gamma) := \sum_{j=1}^\ell \binom{j-1}{i} \eta_j, \text{ for } i = 0, 1, 2.\]
\[ \eta_j(\Gamma) := \eta_j, \qquad \forall j=0, 1, \dots, \ell.\]
\end{notn}

\begin{rem} Note that unlike the case of $P_{\bbeta}(G)$, it is possible
for $\Phi_{\bbeta}(G)$ to be
nonzero even when $G$ is not $\bbeta$-allowable. The reason is that it may
be possible to decompose $G$ into smaller graphs each of which is
$\bbeta$-allowable. Thus, when calculating $\Phi_{\bbeta}(G)$ we cannot 
restrict our attention to the range in which $G$ is $\bbeta$-allowable.
\end{rem}

\section{Deviation from linearity}\label{sec:deviation}

In this section, we analyze the discrepancy between the term 
$Q_{\bbeta}^{\delta}$ which is of primary interest, and the ``linearized''
version $Q_{\delta}(\bbeta)$. Our main result in this direction is
Corollary \ref{cor:qdiff} below. We use the following notation.

\begin{notn}\label{notn:bbetaDelta} 
Given an $h$-transverse polygon $\Delta$, let
$\bbeta(\Delta)$ be the integer sequence with entries consisting of the
widths of $\Delta$ at integer heights, from top to bottom.
\end{notn}

Thus, if $\Delta$ has height $M$, then $\bbeta(\Delta)$ has $M+1$ entries,
which we typically denote by $(\beta_0,\dots,\beta_M)$. If $(\bl_0,\br_0)$
denotes the sequences of edge directions of $\Delta$, with no reordering,
then $\bbeta(\Delta)=\bbeta(d^t,\br_0-\bl_0)$.

In this notation, Corollary \ref{cor:qdiff} is used both in the computation
of $Q_{\bbeta(\Delta)}^{\delta}$ (i.e., the ``unordered case''), and in
the comparison between $Q_{\bbeta}^{\delta}$ and $Q_{\bbeta(\Delta)}^{\delta}$
(i.e., the computation of how reordering affects $Q_{\bbeta}^{\delta}$).

We begin by bounding $\olam_i(\Gamma)$ in terms of the length and cogenus
of $\Gamma$; although for many choices of $\Gamma$ our bounds are far from
sharp, they are sharp when $\Gamma$ is allowed to be arbitrary, and they 
are precisely what we will need in later arguments.

\begin{lem}\label{lem:template-lambda}
Suppose $\Gamma$ is a template of cogenus $\delta$ and length $\ell$. Then
for any $i=1,\dots,\ell$, we have
$$\olam_i \leq \min\{\delta, \delta-(\ell-i)+\epsilon_1(\Gamma), 
\delta+1-i+\epsilon_0(\Gamma)\}.$$
\end{lem}

\begin{proof} We in fact begin by proving that 
$\olam_i(G) \leq \delta(G)$ for any long-edge graph $G$.
Let $m=\lambda_i(G)-\olam_i(G)$, so that $m$ is 
the number of edges from $i-1$ to $i$ in $G$. Let $\lambda'_i$ be the sum
of the multiplicities of these edges. These edges contribute 
$\lambda'_i-m$ to $\delta(G)$. On the other hand, if an edge $e$ has 
length at least $2$ and multiplicity $\rho$, then it contributes at least
$2\rho-1 \geq \rho$ to $\delta(G)$, so the edges not going from $i-1$ to $i$
contribute at least $\lambda_i-\lambda'_i$ to $\delta(G)$.
We thus find 
$$\delta(G) \geq (\lambda'_i-m)+(\lambda_i-\lambda'_i)
=\olam_i,$$
as desired.

To treat the remaining inequalities, we first observe that if $\ell=1$, 
the inequality $\olam_1 \leq \delta$ is the strongest one,
and we have already proved this. We
may thus assume that $\ell>1$, and we first prove that if $i<\ell$, then
we have $\olam_i \leq \delta-(\ell-i)+1$, with strict inequality if 
$\epsilon_1(\Gamma)=0$.
Let $e$ be an edge from $j$ to $\ell'$ chosen so that $\ell'$ is maximal 
subject to the condition that $j \leq i-1$; further suppose that among
the possible choices for $e$ ending at the maximal value of $\ell'$, we
have chosen $e$ so that its weight $\rho$ is maximal. 
Then because $\Gamma$ is a template and $i<\ell$, we have $\ell' > i$.
Let $G$ be the long-edge graph consisting of all the edges of $\Gamma$
other than $e$ which go from $i'$ to $j'$ for some
$i' \leq i-1$ and $j' \geq i$.
Then $\olam_i(G)=\olam_i-\rho$, 
so we conclude from our first inequality that
$\delta(G) \geq \olam_i-\rho$.
On the other hand, $e$ contributes $\rho(\ell'-j)-1$ to $\delta$, and
in order to have length $\ell$, we see by inducting on $\ell-\ell'$ that 
the template condition requires that, in order to cover the interval 
between $\ell'$ and $\ell$, for some $n \geq 0$ there
exist $n$ additional edges of total length at least $\ell-\ell'+n$.
These contribute at least $\ell-\ell'$ to $\delta$, so we conclude that
\begin{align*} \delta 
& \geq \olam_i-\rho+\rho(\ell'-j)-1+\ell-\ell' \\
& = \olam_i+\rho(\ell'-1-j)-1+\ell-\ell' \\
& \geq \olam_i + \ell'-1-j-1+\ell-\ell' \\
& \geq \olam_i+ \ell -i-1,
\end{align*}
as desired. Now, if $\epsilon_1(\Gamma)=0$, this means that an edge ending
at $\ell$ has weight at least $2$. If this edge is not one of the edges
considered in the above inequality, then we certainly obtain strict 
inequality, so assume that it is. If $\ell'<\ell$, this implies that one of 
the $n$ edges considered above has weight at least $2$, so these edges 
contribute strictly greater than $\ell-\ell'$ to
$\delta$. On the other hand, if $\ell'=\ell$, then by the maximality of
$\rho$ we see that $\rho$ is at least as large as any edge of $G$ ending at
$\ell$, so we must have $\rho>1$,
and since $\ell'-1-j\geq\ell'-i>0$, we likewise obtain strict inequality,
as desired.
Similarly, we see that
$\olam_i \leq \delta+2-i$ whenever $i>1$, with strict inequality when
$\epsilon_0(\Gamma)=0$.
The lemma follows, since when $i=\ell$ we have 
$\delta+1-i+\epsilon_0(\Gamma) \leq \delta-(\ell-i)+\epsilon_1(\Gamma)$, 
and when $i=1$ we have
$\delta-(\ell-i)+\epsilon_1(\Gamma) \leq \delta+1-i+\epsilon_0(\Gamma)$.
\end{proof}

In order to illustrate the utility of Lemma \ref{lem:template-lambda},
we also need to briefly examine inequalities on the widths of $\Delta$,
for which we introduce some additional terminology.

\begin{defn} A vertex of $\Delta$ is \textbf{internal} if it is not at
the top or bottom of $\Delta$. An edge is \textbf{internal} if it connects
two internal vertices. An edge is \textbf{extremal} if it connects an
internal vertex to a non-internal vertex.
\end{defn}

Thus, not every edge is internal or extremal, but the only case in which
an edge is neither is when it is at the top or bottom of $\Delta$, or 
connects the top to the bottom.

\begin{prop}\label{prop:allowable-unordered} Let $\Delta$ be an 
$h$-transverse polygon of height $M$ such that 
$\bbeta(\Delta)=(\beta_0,\dots,\beta_{M})$ is
nonconstant, and for $E>0$, suppose that the length of all 
extremal edges of $\Delta$ is at least $E$.
Then we have
$$\beta_i \geq \min\{d^t+i,E,d^b+M-i\}.$$
In particular, we have $\beta_i \geq E$ for $E \leq i \leq M-E$, and
if $\Delta$ has height at least $2E$ (in particular, if $\Delta$ has 
internal vertices), then 
$\beta_i \geq i$ for $i \leq E$ and $\beta_i \geq M-i$ for $i \geq M-E$.
\end{prop}

\begin{proof} The desired inequalities are more or less trivial:
if $\bbeta(\Delta)$ is initially nonincreasing, then it is nonincreasing
everywhere, and must further be strictly 
decreasing for indices greater than or equal to $M-E$, so in this case we 
see in fact that $\beta_i \geq d^b+M-i$ for $i \geq M-E$, and then that
$\beta_i \geq E$ for $i \leq M-E$. Similarly, 
if $\bbeta(\Delta)$ is nondecreasing at the end, 
we check that $\beta_i \geq 
d^t+i$ for $i \leq E$, and then that $\beta_i \geq E$ for $i \geq E$.
Finally, if $\bbeta(\Delta)$ is initially strictly increasing and
strictly
decreasing at the end, then by our hypotheses on the extremal edges we
clearly get the desired inequalities. In the case that 
$M \geq 2E$, we have $i \leq d^b+M-i$ for $i \leq E$, so that
$\min\{d^t+i,E,d^b+M-i\}\geq i$ in this range, and
$M-i \leq d^t+i$ for $i \geq M-E$, so that
$\min\{d^t+i,E,d^b+M-i\}\geq M-i$ in this range.
\end{proof}

We apply these results to relate the numbers $Q_{\bbeta}^{\delta}$ to their
linear approximations $Q_{\delta}(\bbeta)$ in two different contexts.
We first make a few preliminary definitions.

\begin{notn}\label{notn:DQ}
Let $\bv = (0, 1, 2, \dots, \delta).$ For any positive integer $p,$ we 
define 
\[ \DQ(p, \delta) = Q_{p \bv}^{\delta} - Q_\delta(p \bv) 
= \sum_{\Gamma} \mu(\Gamma) \sum_{k=1-\epsilon_0(\Gamma)}^{\delta-\ell(\Gamma) 
+ \epsilon_1(\Gamma)} \Phi_{p \bv} \left( \Gamma_{(k)} \right) 
- \Phi\left( \Gamma_{(k)}, p\bv  \right).\]
where the first summation is over all templates $\Gamma$ of cogenus 
$\delta.$ Let $\DQ(0, \delta) =0.$
\end{notn}

\begin{notn}\label{notn:tdet}
If $d^t=0$, so that $\Delta$ has a unique top vertex $v$, set 
$\tdet(\Delta)= \det(v)$; otherwise, set $\tdet(\Delta)=0$. 
Similarly, if $d^b=0$ and $v$ is the bottom vertex, set 
$\bdet(\Delta)=\det(v)$, and otherwise set $\bdet(\Delta)=0$.
\end{notn}

\begin{defn}\label{defn:conjugate} 
Given a template $\Gamma$, the \textbf{conjugate}
$\nBar{\Gamma}$ of $\Gamma$ is obtained as the image of $\Gamma$ under the 
automorphism of $\ZZ$ given by $n \mapsto \ell(\Gamma)-n$.
\end{defn}

The following lemma then follows directly from the definitions.

\begin{lem}\label{lem:conjugate}
Given a template $\Gamma$, we have
\begin{equation*} 
\ell(\nBar{\Gamma})=\ell(\Gamma), \quad \delta(\nBar{\Gamma})=\delta(\Gamma),
\quad \epsilon_i(\nBar{\Gamma}) = \epsilon_{1-i}(\Gamma) \text{ for } i=0,1, 
\end{equation*}
\begin{equation*}
\eta_0(\nBar{\Gamma})=\eta_0(\Gamma), \quad \eta_i(\nBar{\Gamma}) =
\eta_{\ell(\Gamma)+1-i}(\Gamma) \text{ for } i=1,\dots\ell(\Gamma),
\end{equation*}
\begin{equation*}
\zeta^0(\nBar{\Gamma})=\zeta^0(\Gamma), \quad 
\zeta^1(\Gamma)+\zeta^1(\nBar{\Gamma})=(\ell(\Gamma)-1)\zeta^0(\Gamma).
\end{equation*}
\end{lem}

\begin{cor}\label{cor:qdiff} Suppose we have an $h$-transverse polygon
$\Delta$ such that every extremal edge of $\Delta$ has length at least 
$\delta$, and write $\bbeta(\Delta)=(\beta_0(\Delta),\dots,\beta_M(\Delta)$. 
Then:
\begin{ilist}
\itm if $\bbeta(\Delta)$ is nonconstant and we are given 
$\bbeta=(\beta_0,\dots,\beta_M)$
such that $\beta_i=\beta_i(\Delta)$
for all $i$ with either $i \leq \delta$ or $i\geq M-\delta$,
and $\beta_i \geq \delta$ for $\delta < i < M-\delta$, then
$$Q_{\bbeta(\Delta)}^{\delta} - Q_\delta(\bbeta(\Delta)) =
Q_{\bbeta}^{\delta} - Q_\delta(\bbeta);$$
\itm on the other hand, if further $\Delta$ has height at least $\delta$, 
and each of $d^t,d^b$ is either $0$ or at least $\delta$, then
\begin{align*}
Q_{\bbeta(\Delta)}^{\delta} 
=& Q_\delta(\bbeta(\Delta)) + \DQ(\tdet(\Delta), \delta) 
+\DQ(\bdet(\Delta), \delta).
\end{align*}
\end{ilist}
\end{cor}

\begin{proof} Our main claim is as follows: under the hypotheses of (i),
suppose we are also given a template $\Gamma$ of length $\ell$ and cogenus
$\delta$. Then for all $k$ such that 
$\delta+\epsilon_1(\Gamma) \leq \ell-1+k$ and
$k \leq M-\delta-\epsilon_0(\Gamma)$,
we have that $\Gamma_{(k)}$ is both $\bbeta$-semiallowable and 
$\bbeta(\Delta)$-semiallowable. 
If further $d^t \geq \delta$ (respectively, $d^b \geq \delta$), then
$\Gamma_{(k)}$ is $\bbeta$-semiallowable and $\bbeta(\Delta)$-semiallowable
whenever $0 \leq k \leq M-\delta-\epsilon_0(\Gamma)$ (respectively,
$\delta+\epsilon_1(\Gamma) \leq \ell-1+k \leq M$).

We first observe that by Proposition \ref{prop:allowable-unordered},
we have that $\bbeta(\Delta)$ itself satisfies the hypotheses for 
$\bbeta$, so it is enough to check $\bbeta$-semiallowability.
We then need to show that 
$\olam_i \leq \beta_{i-1+k}$ for $i=1,\dots,\ell$, which we
deduce from Lemma \ref{lem:template-lambda}, using our inequalities on
$\bbeta(\Delta)$ applied in the case $E=\delta$. If we have
$\beta_{i-1+k} \geq \delta$, there is no problem, so it is enough to
consider the cases $i-1+k \leq \delta$ and $i-1+k \geq M-\delta$, for which 
$\beta_{i-1+k}=\beta_{i-1+k}(\Delta)$; we further assume 
$\beta_{i-1+k}<\delta$. Again applying Proposition 
\ref{prop:allowable-unordered},
if $\beta_{i-1+k}=\beta_{i-1+k}(\Delta) \geq d^t+i+k-1 \geq i+k-1$, then 
putting together $\olam_i \leq \delta-(\ell-i)+\epsilon_1(\Gamma)$ with
$\ell-1+k \geq \delta+\epsilon_1(\Gamma)$ gives the desired statement.
The last case is that 
$\beta_{i-1+k}=\beta_{i-1+k}(\Delta) \geq d^b+M+1-i-k \geq M+1-i-k$, from 
which the statement follows similarly using
$\olam_i \leq \delta+1-i+\epsilon_0(\Gamma)$ and 
$k \leq M-\delta-\epsilon_0(\Gamma)$.
Finally, we see that if $d^t \geq \delta$, then the case 
$\beta_{i-1+k} \geq d^t+i+k-1$ is subsumed by $\beta_{i-1+k} \geq \delta$,
so the inequality $\ell-1+k \geq \delta+\epsilon_1(\Gamma)$ is unnecessary,
and if $d^b \geq \delta$, then similarly the inequality
$k \leq M-\delta-\epsilon_0(\Gamma)$ is unnecessary. We thus conclude the
claim.

The first consequence of the claim is that Theorem \ref{thm:linear0} gives
us
\begin{align}\label{eq:QDiff}
Q_{\bbeta(\Delta)}^{\delta} =& Q_\delta(\bbeta(\Delta)) + 
\sum_{\Gamma} \mu(\Gamma) \sum_{k \in I} 
\Phi_{\bbeta(\Delta)} \left( \Gamma_{(k)} \right)-
\Phi\left(\Gamma_{(k)},\bbeta(\Delta)\right), 
\end{align}
and similarly with $\bbeta$ in place of $\bbeta(\Delta)$,
where 
$$I=[1-\epsilon_0(\Gamma),\delta-\ell(\Gamma)+\epsilon_1(\Gamma)] \cup 
[M-\delta+1-\epsilon_0(\Gamma),M-\ell(\Gamma)+\epsilon_1(\Gamma)],$$
and $\Gamma$ ranges over templates of cogenus $\delta$.

We then conclude (i), since the hypothesis that $\beta_i=\beta_i(\Delta)$
for $i \leq \delta$ or $i\geq M-\delta$ implies by Remark 
\ref{rem:linear-support} that
$\Phi_{\bbeta(\Delta)} \left( \Gamma_{(k)} \right)=
\Phi_{\bbeta} \left( \Gamma_{(k)} \right)$ and
$\Phi\left(\Gamma_{(k)},\bbeta(\Delta)\right) =
\Phi\left(\Gamma_{(k)},\bbeta\right)$ for $k \in I$. 

For (ii), if $\bbeta(\Delta)$ is constant, then $d^t \neq 0$, so the width 
of $\Delta$ is at least $\delta$,
and then the desired statement is immediate from Lemma 
\ref{lem:template-lambda} and Theorem \ref{thm:linear0}. We now address
the case that $\bbeta(\Delta)$ is nonconstant. First note that for any 
template $\Gamma$ of cogenus $\delta$, if $d^b=0$ then Lemma 
\ref{lem:conjugate} gives us
$$\sum_{k = M-\delta+1-\epsilon_0(\Gamma)}^{M-\ell(\Gamma)+\epsilon_1(\Gamma)} 
\Phi_{\bbeta(\Delta)} \left( \Gamma_{(k)} \right)-
\Phi\left(\Gamma_{(k)},\bbeta(\Delta)\right)=
\sum_{k = 1-\epsilon_0(\nBar{\Gamma})}^{\delta-\ell(\nBar{\Gamma})+
\epsilon_1(\nBar{\Gamma})} 
\Phi_{p \bv} \left( \nBar{\Gamma}_{(k)} \right)-
\Phi\left(\nBar{\Gamma}_{(k)},p\bv\right),$$
where $\bv$ is as in Notation \ref{notn:DQ} and $p=\bdet(\Delta)$.
We then apply \eqref{eq:QDiff},
observing that the two intervals defining $I$ are disjoint if
the height of $\Delta$ is at least $2\delta$, which is necessarily
the case when $d^t=d^b=0$,
or equivalently when $\tdet(\Delta)$ and $\bdet(\Delta)$ are both nonzero.
The desired formula for this case follows. 
If $d^t=0$ but $d^b \geq \delta$, then the above claim gives us 
semiallowability for
$\delta+\epsilon_1(\Gamma) \leq \ell(\Gamma)-1+k \leq M$, so we can 
replace $I$ in \eqref{eq:QDiff} by 
$[1-\epsilon_0(\Gamma),\delta-\ell(\Gamma)+\epsilon_1(\Gamma)]$, and noting
that $\DQ(\bdet(\Delta), \delta)=0$ in this case by convention, we 
again get the desired formula. The cases that $d^t \geq \delta$ but
$d^b=0$ and that $d^t,d^b \geq \delta$ are similar, so we conclude that
the stated formula holds in all cases.
\end{proof}

\section{Formulas for the linearization}\label{sec:no-reorder}

In this section, we analyze $Q_{\delta}(\bbeta)$ in detail,
giving an explicit combinatorial formula for it which we then apply
to give a formula for $Q_{\bbeta(\Delta)}^{\delta}$. 

We begin with the definitions of the invariants which 
form the basis for our main formulas.

\begin{notn}\label{notn:coefs}
Given $\delta \geq 1$, set
\begin{align*}
  A(\delta) =& \frac{1}{2} \sum \mu(\Gamma) \zeta^0(\Gamma), \\
  L(\delta) :=& -\frac{1}{2} \sum \mu(\Gamma) \zeta^0(\Gamma)(\ell(\Gamma) - \epsilon_0(\Gamma)-\epsilon_1(\Gamma)), \\
  H(\delta) := & \sum \mu(\Gamma) \left( \eta_0(\Gamma) + \zeta^0(\Gamma)(\ell(\Gamma) - \epsilon_0(\Gamma)-\epsilon_1(\Gamma)) \right), \\
  D(\delta) := &- \sum \mu(\Gamma) \left( \zeta^2(\Gamma) + \zeta^1(\Gamma)(1-\epsilon_0(\Gamma) ) \right), \\
  C(\delta) := & - \sum \mu(\Gamma) \eta_0(\Gamma) (\ell(\Gamma) - \epsilon_0(\Gamma)-\epsilon_1(\Gamma)), 
\end{align*}
  where all the summations are over all the templates of cogenus $\delta.$
\end{notn}

We use $\zeta^1$ and $\eta_0$ to compute a special case of $\DQ(p,\delta)$,
as follows.

\begin{lem}\label{lem:DQ-pgedelta}
For any $p \ge \delta,$
\[\DQ(p, \delta) = - \sum_{\Gamma: \epsilon_0(\Gamma)=1}\mu(\Gamma) \left(p\zeta^1(\Gamma) + \eta_0(\Gamma)\right).\]
\end{lem}

\begin{proof}
With $\bv$ as in Notation \ref{notn:DQ}, because $p \geq \delta$, we have 
that $\Gamma_{(k)}$ is $p\bv$-semiallowable for $k \geq 1$.
On the other hand, if $k=0$, we observe that 
$\Phi_{p\bv}(\Gamma_{(k)}) = \Phi_{p\bv}(\Gamma) = 0$, since 
any partition of $\Gamma$ must contain at least one graph having
an edge incident to $0$. Then by Theorem \ref{thm:linear0}, we find
$$\DQ(p, \delta) =
-\sum_{\Gamma: \epsilon_0(\Gamma)=1} \mu(\Gamma) \Phi(\Gamma, p\bv),$$
which yields the desired formula.
\end{proof}

We next make the following definitions, which are motivated by
sequences of the form $\bbeta(\Delta)$, but make sense (and are important)
more generally.

\begin{notn}\label{notn:areaetc}
  Suppose $\bbeta= (\beta_0, \dots, \beta_M) = \bbeta(\bd),$ where $\bd=(d_0, d_1, \dots, d_M).$ Define
  \begin{align*}
    \area(\bbeta) =&  \beta_0 + 2\beta_1 + \cdots + 2 \beta_{M-1} + \beta_M \\
    \LL(\bbeta) =& \beta_0 + \beta_M + 2M \\
    \height(\bbeta) =& M \\
    \idet(\bbeta) =& d_1 - d_M \\
  \end{align*}
\end{notn}

The main formula of this section is then as follows.

\begin{prop}\label{prop:QlinwH}
Suppose $\bbeta=(\beta_0, \beta_1, \dots, \beta_M)= \bbeta(\bd),$ where $\bd=(d_0, d_1, \dots, d_M)$ satisfying   
  \[ d_1 = d_2 = \cdots = d_{\delta-1} \text{ and } d_M = d_{M-1} = \cdots = d_{M+2-\delta}.\]
  Then
\[ Q_\delta(\bbeta) = A(\delta) \cdot \area(\bbeta) + L(\delta) \cdot \LL(\bbeta) + H(\delta) \cdot \height(\bbeta) + D(\delta) \cdot \idet(\bbeta) + C(\delta).\]
\end{prop}

Implicit in the hypotheses of Proposition \ref{prop:QlinwH} is that
$M \geq \delta-1$; we also have $M \geq 1$ for any polygon $\Delta$.
The proof is largely straightforward, 
with the main trick being to consider each template $\Gamma$ together 
with its conjugate. A preliminary lemma is the following.

\begin{lem}\label{lem:sumPhi}
Suppose $\Gamma$ is a template of length $\ell$
 and let $a,b$ be integers satisfying $0 \le a \le b+1 \le M-\ell+2.$
Suppose $\bbeta=(\beta_0, \beta_1, \dots, \beta_M)= \bbeta(\bd),$ where 
$\bd=(d_0, d_1, \dots, d_M)$ satisfies 
\[ d_{a+1} = d_{a+2} = \cdots = d_{a+\ell-2} = p, \quad d_{b+\ell-1} = d_{b+\ell-2} = \cdots = d_{b+2} = q.\]
Then
\begin{align*}
  \sum_{k=a}^b \Phi\left(  \Gamma_{(k)} , \bbeta \right)  
=& \ \eta_0(\Gamma) (b-a+1) 
+ \zeta^0(\Gamma) \left(\sum_{k=a}^{b+\ell-1} \beta_k\right) \\
& - \zeta^1(\Gamma) \beta_a - \zeta^1(\nBar{\Gamma}) \beta_{b+\ell-1} 
  - \zeta^2(\Gamma) p + \zeta^2(\nBar{\Gamma}) q. 
 \end{align*}
 \end{lem}

\begin{proof}
 It follows from Remark \ref{rem:linear-support} that 
\begin{align*}
 \Phi\left(\Gamma_{(k)}, \bbeta\right) =&  \Phi\left(\Gamma, (\beta_k, \beta_{k+1}, \dots, \beta_{k+\ell-1})\right) \\
 =& \eta_0(\Gamma) + \eta_1(\Gamma) \beta_k + \eta_2(\Gamma) \beta_{k+1} + \cdots + \eta_\ell(\Gamma) \beta_{k+\ell-1},
 \end{align*}
and then
it follows from the hypotheses that for any $1 \le k \le \ell-1,$
\[ \beta_{a+k-1} = \beta_a + (k-1) p, \quad \text{ and } \quad 
\beta_{b+\ell-k} = \beta_{b+\ell-1} - (k-1) q.\]
Thus, it follows from the above together with Lemma \ref{lem:conjugate} that
 \begin{align*}
 \sum_{k=a}^b \Phi\left(  \Gamma_{(k)} , \bbeta \right)  =& \ \sum_{k=a}^b\left(\eta_0(\Gamma) + \eta_1(\Gamma) \beta_k + \eta_2(\Gamma) \beta_{k+1} + \cdots + \eta_\ell(\Gamma) \beta_{k+\ell-1}  \right) \\ 
 =& \eta_0(\Gamma) (b-a+1) 
+ \zeta^0(\Gamma) \left(\sum_{k=a}^{b+\ell-1} \beta_k\right) \\
 & - \sum_{k=1}^{\ell-1} \left(\sum_{i=k+1}^{\ell} \eta_{i}(\Gamma)\right) \cdot \beta_{a+k-1} - \sum_{k=1}^{\ell-1} \left(\sum_{i=k+1}^{\ell} \eta_{\ell+1-i}(\Gamma)\right) \cdot \beta_{b+\ell-k} \\
  =& \ \eta_0(\Gamma) (b-a+1) 
+ \zeta^0(\Gamma) \left(\sum_{k=a}^{b+\ell-1} \beta_k\right) 
- \zeta^1(\Gamma) \beta_a - \zeta^1(\nBar{\Gamma}) \beta_{b+\ell-1} \\
  & - p \sum_{k=1}^{\ell-1} (k-1) \left(\sum_{i=k+1}^\ell \eta_{i}(\Gamma)\right) + q \sum_{k=1}^{\ell-1} (k-1) \left(\sum_{i=k+1}^{\ell} \eta_{i}(\nBar{\Gamma})\right),
 \end{align*}
and the desired formula follows.
\end{proof}

We can now prove the asserted formula.

\begin{proof}[Proof of Proposition \ref{prop:QlinwH}]
  Observing that
  \[ Q_\delta(\bbeta) = \frac{1}{2}\sum_{\Gamma} \left( Q_\Gamma(\bbeta) + Q_{\nBar{\Gamma}}(\bbeta) \right) ,\]
  where the summation is over all templates $\Gamma$ of cogenus $\delta$,
we see that it is enough to prove that for each such $\Gamma$, we have

 \begin{align*} 
\frac{1}{2\mu(\Gamma)}( Q_\Gamma(\bbeta) +& Q_{\nBar{\Gamma}}(\bbeta)) \\
=& - \eta_0(\Gamma)(\ell(\Gamma)-\epsilon_0(\Gamma)-\epsilon_1(\Gamma)) + \frac{1}{2} \zeta^0(\Gamma) \cdot \area(\bbeta)  \\
 & - \frac{1}{2} \zeta^0(\Gamma)  (\ell(\Gamma)-\epsilon_0(\Gamma)-\epsilon_1(\Gamma)) \cdot \LL(\bbeta) \\
 & + (\eta_0(\Gamma) +\zeta^0(\Gamma)(\ell(\Gamma)-\epsilon_0(\Gamma)-\epsilon_1(\Gamma) ) \cdot \height(\bbeta) \\
 & - \frac{1}{2}  \left(\zeta^2(\Gamma) + \zeta^2(\nBar{\Gamma}) + \zeta^1(\Gamma)(1-\epsilon_0(\Gamma)) +\zeta^1(\nBar{\Gamma})(1-\epsilon_0(\nBar{\Gamma})) \right) \cdot \idet(\bbeta).
\end{align*}

Now, for any $\Gamma$ (and in particular also for $\nBar{\Gamma}$) we have
$1-\epsilon_0(\Gamma) \leq M-\ell(\Gamma)+\epsilon_1(\Gamma)+1$, and
 \[ (1-\epsilon_0(\Gamma)) + \ell(\Gamma) -2 \le \delta-1, \quad\text{ and } \quad (M-\ell(\Gamma)+\epsilon_1(\Gamma))+2 \ge M+2-\delta,\]
so applying Lemma \ref{lem:sumPhi} with $a=1-\epsilon_0(\Gamma)$ and
$b=M-\ell(\Gamma)+\epsilon_1(\Gamma)$ together with Lemma \ref{lem:conjugate},
and recalling that $d_1=\beta_1-\beta_0$ and $d_M=\beta_M-\beta_{M-1}$,
we get
  \begin{align*}
    & \frac{1}{2\mu(\Gamma)}\left( Q_\Gamma(\bbeta) + Q_{\nBar{\Gamma}}(\bbeta) \right) = \frac{1}{2} \left( \sum_{k=1-\epsilon_0(\Gamma)}^{M-\ell(\Gamma) + \epsilon_1(\Gamma)} \Phi\left( \Gamma_{(k)}, \bbeta \right) + \sum_{k=1-\epsilon_0(\nBar{\Gamma})}^{M-\ell(\nBar{\Gamma}) + \epsilon_1(\nBar{\Gamma})} \Phi\left( \nBar{\Gamma}_{(k)}, \bbeta \right)  \right) \\
     =& \eta_0(\Gamma)(M-\ell(\Gamma)+\epsilon_0(\Gamma)+\epsilon_1(\Gamma)) + \frac{1}{2} \zeta^0(\Gamma) \left( \sum_{k=1-\epsilon_0(\Gamma)}^{M + \epsilon_1(\Gamma)-1} \beta_k +  \sum_{k=1-\epsilon_1(\Gamma)}^{M+ \epsilon_0(\Gamma)-1} \beta_k \right) \\
     & - \frac{1}{2} \left( \zeta^1(\Gamma) (\beta_{1-\epsilon_0(\Gamma)} + \beta_{M+\epsilon_0(\Gamma)-1}) +  \zeta^1(\nBar{\Gamma}) (\beta_{1-\epsilon_1(\Gamma)} + \beta_{M+\epsilon_1(\Gamma)-1} ) \right) \\
     & - \frac{1}{2}(d_1-d_M) (\zeta^2(\Gamma) + \zeta^2(\nBar{\Gamma})).
  \end{align*}
  We see that
\[ \sum_{k=1-\epsilon_0(\Gamma)}^{M + \epsilon_1(\Gamma)-1} \beta_k +  \sum_{k=1-\epsilon_1(\Gamma)}^{M+ \epsilon_0(\Gamma)-1} \beta_k
= \area(\bbeta) +(\epsilon_0(\Gamma)+\epsilon_1(\Gamma)-1)(\beta_0+\beta_M),\]
and
\begin{align*}
& \zeta^1(\Gamma) (\beta_{1-\epsilon_0(\Gamma)} + \beta_{M+\epsilon_0(\Gamma)-1}) +  \zeta^1(\nBar{\Gamma}) (\beta_{1-\epsilon_1(\Gamma)} + \beta_{M+\epsilon_1(\Gamma)-1} ) \\
=& \zeta^1(\Gamma)(\beta_0+\beta_M + (d_1-d_M)(1-\epsilon_0(\Gamma))) +  \zeta^1(\nBar{\Gamma})(\beta_0+\beta_M + (d_1-d_M)(1-\epsilon_1({\Gamma}))) \\
=& (\ell(\Gamma)-1) \zeta^0(\Gamma)(\beta_0+\beta_M) + ( \zeta^1(\Gamma)(1-\epsilon_0(\Gamma)) +\zeta^1(\nBar{\Gamma})(1-\epsilon_0(\nBar{\Gamma}))) (d_1-d_M),
\end{align*}
and the desired formula follows by rearranging terms.
\end{proof}

Finally, putting together Proposition \ref{prop:QlinwH} with
Corollary \ref{cor:qdiff} (ii), we conclude the following.

\begin{cor}\label{cor:wH}
Suppose that all extremal edges of $\Delta$ have length at least $\delta$,
that $\Delta$ has height at least $\delta$, and that each of $d^t$ and
$d^b$ are at least $\delta$ if they are nonzero.
  Then
\begin{align*}
Q_{\bbeta(\Delta)}^{\delta} 
=& A(\delta) \cdot \area(\Delta) + L(\delta) \cdot \LL(\Delta)  
+ H(\delta) \cdot \height(\Delta) + D(\delta) \cdot \idet(\Delta)  \\
&+ C(\delta) + \DQ(\tdet(\Delta), \delta) +\DQ(\bdet(\Delta), \delta),
\end{align*}
where $\area(\Delta)$ is the normalized area of $\Delta$ (that is, twice
the area), $\LL(\Delta)$ is the lattice length of $\Delta$, $\height(\Delta)$
is the height of $\Delta$, and $\idet(\Delta)$ is the sum of the determinants
of the internal vertices of $\Delta$.
\end{cor}

Note that $\area(\Delta)=\area(\bbeta(\Delta))$, 
$\LL(\Delta)=\LL(\bbeta(\Delta))$, $\height(\Delta)=\height(\bbeta(\Delta))$,
and $\idet(\Delta)=\idet(\bbeta(\Delta))$.

\section{The behavior of reorderings}\label{sec:reorder}

The purpose of this section is to address how to remove the summing over 
reorderings
$(\bl,\br)$ in the definition of $N^{\Delta,\delta}$, giving formulas
for both $\cN$ and $\cQ$ stated solely in terms of the default ordering
$(\bl_0,\br_0)$ corresponding to the unpermuted edge directions. 

The main results of this section are then as follows. We begin by
analyzing the effect of reordering $(\bl,\br)$ on $Q_{\bbeta}^{\delta}$.

\begin{lem}\label{lem:reversals} 
If every extremal edge of $\Delta$ has length at least $\delta(\bl,\br)$,
then $\bbeta(d^t,\br-\bl) \in \ZZ_{\geq 0}^{M+1}$. If further every
extremal edge has length at least $\delta+\delta(\bl,\br)$, then
\[ Q_{\bbeta(d^t, \br - \bl)}^{\delta} = 
Q_{\bbeta(\Delta)}^{\delta} - 2A(\delta) \cdot \delta(\bl, \br).\]
\end{lem}

Because it is $\cN^{\Delta}(t)$ rather than $\cQ^{\Delta}(t)$ which is defined 
in terms of summing over all choices of $(\bl,\br)$, some additional work
is necessary in order to extract the most useful formulas from Lemma 
\ref{lem:reversals}. 

\begin{notn}\label{notn:AP}
Write $\cP(x)=\sum_{n \ge 0} p(n) x^n,$ the generating function of the 
partitions of numbers. Set
$$\cA(t)=\exp\left(-\sum_{\delta \geq 1} 2A(\delta) t^{\delta}\right),$$
\end{notn}

Using Lemma \ref{lem:reversals}, we obtain a formula for 
$N^{\Delta,\delta}$ which does not involve summing over reorderings
$(\bl,\br)$.

\begin{notn}\label{notn:ellDelta}
Define $\ell(\Delta)$ to be the minimum of the 
lengths of the internal edges and extremal edges of $\Delta$.
\end{notn}

We make the convention that $\ell(\Delta)=+\infty$ if $\Delta$ has no
internal vertices (and hence no extremal or internal edges).

\begin{cor}\label{cor:reversals-N}
If $\ell(\Delta) \geq \delta$, and we denote by $v'_i$ the number of
internal vertices of $\Delta$ having determinant $i$ then
\[ N^{\Delta, \delta} = 
\left[ t^{\delta} \right] \left(\cN(\bbeta(\Delta),t) \cdot 
\left(\prod_{i=1}^{\delta} \left(\cP\left(\left( t\cA(t) \right)^i
\right)\right)^{v'_i}\right)\right).\]
\end{cor}

The above formula may be of independent interest, but for our purposes
is used only to obtain the below formula for $Q^{\Delta,\delta}$. We
first introduce the following coefficients:

\begin{equation}\label{equ:bdeltai}
b_{\delta,i} 
:= \left[ t^{\delta} \right] \left( \sum_{n \ge1} \frac{\sum_{d |n} d}{n} 
\cdot (t \cA(t))^{i \cdot n} \right).
\end{equation}

We finally conclude the following corollary.

\begin{cor}\label{cor:reversals-Q}
If $\ell(\Delta) \geq \delta$, and we denote by $v'_i$ the number of
internal vertices of $\Delta$ having determinant $i,$ then
\[ Q^{\Delta, \delta} 
= Q_{\bbeta(\Delta)}^{\delta} + 
\sum_{i=1}^{\delta} b_{\delta,i} v'_i.\]
\end{cor}

We now give the proof of the stated results, starting with some simple
observations.

\begin{prop}\label{prop:swap} Suppose that $(\bl',\br')$ is obtained from
$(\bl,\br)$ by swapping a single adjacent reversed pair in $-\bl$ or $\br$, 
which differ from one another by $d$.
Then
$$\delta(\bl',\br')=\delta(\bl,\br)-d.$$
\end{prop}

\begin{defn} Given an internal vertex $v$ on the right
side of $\Delta$, suppose that $N_1<N_2$ are the levels of the vertices
of $\Delta$ immediately above and below $v$. We say a reordering 
$(\bl,\br)$ is \textbf{$v$-local} $\bl$ is equal to the (ordered) set of 
left directions of $\Delta$, and if $\br$ can be obtained from the set of 
right directions of $\Delta$ using permutations supported in $(N_1,N_2]$. 

We make the same definition if $v$ is on the left side of $\Delta$, with 
the roles of left and right switched.
\end{defn}

Here we say a vertex has level $N$ if $N$ is the vertical distance to
the vertex from the top of $\Delta$.

We then also observe the following.

\begin{prop}\label{prop:reorder-decomp} Suppose that every internal edge
of $\Delta$ has length at least $\delta$ for some $\delta \geq 0$.
Then there is a bijection between reorderings $(\bl,\br)$ with
$\delta(\bl,\br)=\delta$ and tuples of $v$-local reorderings 
$((\bl_v,\br_v))_{v \in \Delta}$ with $v$ varying over internal vertices
of $\Delta$ and $\sum_v \delta(\bl_v,\br_v) = \delta$.
\end{prop}

\begin{proof} The basic observation is the following: if $v$ is an 
internal vertex of determinant $d$ and level $N$, suppose without loss
of generality that $v$ is on the right side of $\Delta$, and $r>r'$ are
the adjacent edge directions. Then if $(\bl,\br)$ is a reordering with $r'$
appearing in index $N+1-i$, we must have
\begin{equation}\label{eq:local-ineq}
\delta(\bl,\br) \geq di \geq i.
\end{equation}
Indeed, in this case we must make at least $i$ adjacent swaps involving $r'$
in order to return to the default ordering, so \eqref{eq:local-ineq} follows
from Proposition \ref{prop:swap}. Similarly, if $r$ occurs in index at 
least $N+i$, we see that \eqref{eq:local-ineq} still holds. 

We then conclude that if $v,v'$ are adjacent internal vertices of $\Delta$ 
of levels $N,N'$ respectively, and $(\bl,\br)$ is a reordering with 
$N'-N \geq \delta(\bl,\br)$, then if $v,v'$ are on the right with adjacent
edge directions $r>r'>r''$, we cannot have any $r''$ occurring before any
$r$ in $\br$, and similarly if $v,v'$ are on the left. Thus, 
$(\bl,\br)$ decomposes naturally into a tuple of $v$-local reorderings.
Similarly, if we have a tuple $((\bl_v,\br_v))_v$ of $v$-local reorderings,
and $v,v'$ are as above, we see that if
with $N'-N \geq \delta(\bl_v,\br_v)+\delta(\bl_{v'},\br_{v'})$, then if
$v,v'$ are on the right with adjacent edge directions $r>r'>r''$, 
the smallest index
in which $r''$ occurs in $\br_{v'}$ must be strictly greater than the largest
index in which $r$ occurs in $\br_v$.
The same holds on the left, and applying this to all adjacent pairs of 
internal vertices implies that we can combine the $(\bl_v,\br_v)$ to get a 
single reordering $(\bl,\br)$. The proposition follows.
\end{proof}

We can now prove our main lemma.

\begin{proof}[Proof of Lemma \ref{lem:reversals}]
First note that nontrivial reorderings only occur when $\Delta$
has internal vertices, so we may assume that $\Delta$
has internal vertices, and in particular that $\bbeta(\Delta)$ is
nonconstant. Write $\bbeta(d^t,\br-\bl)=(\beta_0,\dots,\beta_M)$ and
$\bbeta(\Delta)=(\beta_0(\Delta),\dots,\beta_M(\Delta))$.
Next, observe that inductively applying Proposition \ref{prop:swap} to 
adjacent swaps shows that we have
\begin{equation}\label{eq:alt-delta}
\delta(\bl,\br)=\sum_{j=0}^{M} \left(\beta_j(\Delta)-\beta_j\right),
\end{equation}
and moreover, by considering the effect on $\bbeta$ of adjacent swaps,
we also see that each term on the right is nonnegative.

We then claim that, for any $\delta \geq 0$, if every extremal edge of
$\Delta$ has length at least $\delta+\delta(\bl,\br)$, we have
$\beta_i=\beta_i(\Delta)$ when $i\leq \delta$ or 
$i \geq M-\delta$, and 
$\beta_i\geq \delta$ when $\delta < i < M-\delta$.
Applying this claim in the case $\delta=0$ proves the first statement of
the lemma. Next, given $\delta \geq 1$, putting the claim together with 
Corollary \ref{cor:qdiff} (i) reduces the lemma to showing that
\[ Q_{\delta}(\bbeta(d^t, \br - \bl)) = 
Q_{\delta}(\bbeta(\Delta)) - 2A(\delta) \cdot \delta(\bl, \br).\]
We then apply Proposition \ref{prop:QlinwH} to both sides,
noting that we have $\LL(\bbeta(\Delta))=\LL(\bbeta(d^t,\br-\bl))$,
$\height(\bbeta(\Delta))=\height(\bbeta(d^t,\br-\bl))$, and
$\idet(\bbeta(\Delta))=\idet(\bbeta(d^t,\br-\bl))$.
We thus see that we want to show that
$$\area(\bbeta(d^t,\br-\bl))=\area(\bbeta(\Delta))-2\delta(\bl,\br),$$
which follows immediately from \eqref{eq:alt-delta}.
It is thus enough to prove the claim.

Let $N$ be the level of the uppermost internal vertex, so that we have
assumed that $N \geq \delta+\delta(\bl,\br)$.
Suppose first that $\beta_i(\Delta)-\beta_i =n>0$ for some $i \leq N$.
Then because the right edge directions between $0$ and $N$ are maximal
and the left ones are minimal, we see that 
$\beta_{i'}(\Delta)-\beta_{i'} \geq n$ for all $i'=i,\dots,N$,
so \eqref{eq:alt-delta} implies that
\begin{equation}\label{eq:delta-ineq}
\delta(\bl,\br) \geq n (N+1-i).
\end{equation}
In particular, we have $i \geq N+1-\delta(\bl,\br) \geq \delta+1$.
Considering similarly the lowest internal vertex, we see also that
$i < M-\delta$, giving the first statement we wished to show.

Now, applying Proposition \ref{prop:allowable-unordered} in the
case $E=\delta(\bl,\br)+\delta$ gives us that if
$\delta(\bl,\br)+\delta \leq i \leq M-\delta(\bl,\br)-\delta$, then
$\beta_i(\Delta) \geq \delta(\bl,\br)+\delta$, and
\eqref{eq:alt-delta} then shows that $\beta_i \geq \delta$ for these
values of $i$. It thus remains to show that 
when $\delta < i < \delta(\bl,\br)+\delta$
or $M-\delta(\bl,\br)-\delta < i < M-\delta$, we still have
$\beta_{i} \geq \delta$.
First consider
the case $\delta < i < \delta(\bl,\br)+\delta$. In this range, because
$M \geq 2(\delta(\bl,\br)+\delta)$, we have $i \leq M-i$, so
Proposition \ref{prop:allowable-unordered} gives us that
$\beta_{i}(\Delta) \geq i$. Thus, if 
$\beta_{i}(\Delta)-\beta_{i} =n$, we want to prove that
$i-n \geq \delta$. We will in fact show this for 
$\delta<i \leq \delta(\bl,\br)+\delta$.
According to \eqref{eq:delta-ineq}, we have
$n \leq \frac{\delta(\bl,\br)}{N+1-i}$,
so we wish to show that
$\frac{\delta(\bl,\br)}{N+1-i} \leq i-\delta$,
or equivalently,
$$ \delta(\bl,\br)  \leq (N+1-i) (i-\delta).$$ 
Since the righthand side is quadratic in $i$ with negative leading 
coefficient, to check the desired inequality for 
$i=\delta+1,\dots,\delta(\bl,\br)-\delta$,
it is enough to check it at the endpoints. For $i=\delta+1$, we get
$$(N+1-1-\delta)\cdot (\delta+1-\delta)= N-\delta \geq \delta(\bl,\br),$$
as desired. For $i=\delta(\bl,\br)+\delta$, we similarly get
$$(N+1-\delta(\bl,\br)-\delta) (\delta(\bl,\br)+\delta-\delta) =
(N+1-\delta(\bl,\br)-\delta) \delta(\bl,\br) \geq 1 \cdot \delta(\bl,\br),$$
proving the claim for $\delta<i<\delta(\bl,\br)+\delta$. 
The argument for
$M-\delta(\bl,\br)-\delta < i-1+k < M-\delta$ is the same, considering
the level of the lowest internal vertex rather than the highest. We have
thus proved the claim, and the lemma.
\end{proof}

We conclude by giving the proof of our corollaries to Lemma
\ref{lem:reversals}.

\begin{proof}[Proof of Corollaries \ref{cor:reversals-N} and 
\ref{cor:reversals-Q}]
Let $\ell^e(\Delta)$ be the minimum length of the extremal edges of 
$\Delta$, with the convention that $\ell^e(\Delta)=+\infty$ if $\Delta$
has no internal vertices.
We may rephrase Lemma \ref{lem:reversals} as saying that 
$$\cQ(\bbeta(d^t,\br-\bl),t) \equiv
\cQ(\bbeta(\Delta),t)
-\delta(\bl,\br) \sum_{\delta \geq 1} 2A(\delta) t^{\delta}
\pmod{t^{\ell^e(\Delta)+1-\delta(\bl,\br)}}$$
whenever $\delta(\bl,\br)\leq \ell^e(\Delta)$.
Taking exponentials on both sides, we find
\begin{align*} \cN(\bbeta(d^t,\br-\bl),t) 
& \equiv \cN(\bbeta(\Delta),t)
/\exp \left(\delta(\bl,\br) \sum_{\delta \geq 1} 2A(\delta) t^{\delta}\right)
\pmod{t^{\ell^e(\Delta)+1-\delta(\bl,\br)}} \\
& = \cN(\bbeta(\Delta),t) \cdot \cA(t)^{\delta(\bl,\br)}. 
\end{align*}
If we multiply by $t^{\delta(\bl,\br)}$ and 
sum over choices of $(\bl,\br)$ with 
$\bbeta(d^t,\bl-\br) \in \ZZ_{\geq 0}^{M+1}$,
we find
\begin{align*}
\cN^{\Delta}(t) & := \sum_{(\bl,\br):\bbeta(d^t,\br-\bl) \in \ZZ_{\geq 0}^{M+1}}
\cN(\bbeta(d^t,\br-\bl),t)t^{\delta(\bl,\br)} \\
& \equiv \sum_{(\bl,\br):\delta(\bl,\br) \leq \ell^e(\Delta)}
\cN(\bbeta(d^t,\br-\bl),t)t^{\delta(\bl,\br)} 
\pmod{t^{\ell^e(\Delta)+1}}
\\
& \equiv 
\cN(\bbeta(\Delta),t) \cdot 
\left(\sum_{(\bl,\br):\delta(\bl,\br) \leq \ell^e(\Delta)} 
\left(t\cA(t)\right)^{\delta(\bl,\br)} \right)
\pmod{t^{\ell^e(\Delta)+1}},
\end{align*}
where the first congruence is a consequence of the assertion of
Lemma \ref{lem:reversals} that if $\delta(\bl,\br) \leq \ell^e(\Delta)$,
then $\bbeta(d^t,\bl-\br) \in \ZZ_{\geq 0}^{M+1}$.

Then, because $\ell(\Delta) \leq \ell^e(\Delta)$ by definition,
Proposition \ref{prop:reorder-decomp} gives us
$$\cN^{\Delta}(t) \equiv
\cN(\bbeta(\Delta),t) \cdot 
\left(\prod_{v \in \Delta} \sum_{(\bl_v,\br_v)} 
\left(t\cA(t)\right)^{\delta(\bl_v,\br_v)} \right)
\pmod{t^{\ell(\Delta)+1}},$$
where the sum is over $v$-local reorderings $(\bl_v,\br_v)$. Note that
although \textit{a priori} we should be considering tuples $(\bl_v,\br_v)$
of $v$-local reorderings such that 
$\sum_v \delta(\bl_v,\br_v) \leq \ell(\Delta)$, the tuples with
$\sum_v \delta(\bl_v,\br_v) > \ell(\Delta)$ do not contribute to
the righthand side when considered modulo $t^{\ell(\Delta)+1}$.

Now, we claim that
$$\sum_{(\bl_v,\br_v)} x^{\delta(\bl_v,\br_v)}\equiv \cP(x^d) 
\pmod{x^{\ell(\Delta)+1}},$$
where $d=\det(v)$.
To prove this, for any $m \leq \ell(\Delta)$, we construct a bijection 
between $v$-local reorderings
$(\bl_v,\br_v)$ with $\delta(\bl_v,\br_v)=d m$, and partitions of $m$.
We construct this bijection inductively as follows: given $(\bl_v,\br_v)$, 
let $\lambda$ be the empty partition, suppose without loss of generality 
that $v$ is on the right side of $\Delta$, and let $i$ be minimal 
with $r_{i+1}>r_i$. Let 
$i'$ be minimal with $r_{i'}=r_i$, so that $r_{i'}=r_{i'+1}=\dots=r_i$.
Then let $(\bl'_v,\br'_v)$ be obtained from $(\bl_v,\br_v)$ by swapping 
$r_{i'}$ with $r_{i+1}$; this can be accomplished by a sequence of 
$i+1-i'$ adjacent swaps, so by Proposition \ref{prop:swap}, we see that
$\delta(\bl'_v,\br'_v) = d(m-(i+1-i'))$. Prepend $i+1-i'$ to $\lambda$,
and then repeat the process with $(\bl',\br')$ in place of
$(\bl,\br)$, until $m=0$ and $(\bl,\br)=(\bl_0,\br_0)$. It is then clear
that $\lambda$ is a composition of $m$, but since $i$ strictly increases at
each step, while $i'$ only increases by $1$, we see that $\lambda$ is
weakly decreasing, so gives a partition of $m$, as desired. One checks 
easily that this process is invertible when $m \leq \ell(\Delta)$, 
proving the claim. This yields
$$\cN^{\Delta}(t) \equiv
\cN(\bbeta(\Delta),t) \cdot 
\prod_{i \geq 1} \left(\cP \left((t\cA(t))^i \right) \right)^{v'_i}
\pmod{t^{\ell(\Delta)+1}},$$
which, noting that the product on the righthand side only affects the
$t^{\delta}$ term for $i \leq \delta$, is equivalent to the statement 
of Corollary \ref{cor:reversals-N}.

If we take the formal logarithm on both sides we obtain
$$\cQ^{\Delta}(t) \equiv
\cQ(\bbeta(\Delta),t) + 
\sum_{i \geq 1} v'_i \log \cP \left((t\cA(t))^i \right) 
\pmod{t^{\ell(\Delta)+1}},$$
and noting that 
\begin{equation}\label{eq:log-P} 
\log \cP(x) = \log \left(\prod_{i \geq 1} \frac{1}{1-x^i}\right) 
= \sum_{n \ge 1} \frac{\sum_{d | n} d}{n} \cdot x^n\end{equation}
and that once again only terms with $i \leq \delta$ contribute to
the $t^{\delta}$ term,
we also obtain the statement of Corollary \ref{cor:reversals-Q}.
\end{proof}

\section{Proof of the main theorem}\label{sec:geometry}

In this section, we use some tricks involving the geometric origins of
our formulas in order to compute some of the combinatorial terms which
we had derived previously. Specifically, we use that if $\Delta$ is
$h$-transverse, and a $90$-degree rotation $\Delta'$ of $\Delta$ is 
also $h$-transverse, then these correspond to isomorphic polarized
toric varieties, so we must have $\cN^{\Delta}(t)=\cN^{\Delta'}(t)$, and
correspondingly $\cQ^{\Delta}(t)=\cQ^{\Delta'}(t)$. We then replace our
variables with algebrogeometric ones compatible with Tzeng's formula,
and complete the proof of our main theorems.

\begin{lem}\label{lem:h-0}
For any $\delta,$ we have $H(\delta) = 0.$
\end{lem}

\begin{proof}
Consider $\Delta$ to be the rectangle with horizontal edges of length $a$ 
and vertical edges of length $b,$ where $a \neq b$ and $a, b \ge \delta.$ 
Corollary \ref{cor:wH} says that
$[t^\delta] \cQ(\bbeta(\Delta),t)$ is given by
\begin{multline*}
A(\delta) \cdot \area(\Delta) + L(\delta) \cdot \LL(\Delta) 
+ H(\delta) \cdot \height(\Delta) + D(\delta) \cdot \idet(\Delta) 
+ C(\delta) \\
= A(\delta) \cdot (2ab) + L(\delta) \cdot (2a+2b) + H(\delta) \cdot b 
+ C(\delta),
\end{multline*}
but since there are no nontrivial reorderings of left or right directions, 
we have
\[ \cQ(\bbeta(\Delta),t) = \log(\cN(\bbeta(\Delta),t)) 
= \log(\cN^{\Delta}(t)). \]
Thus, if we obtain $\Delta'$ by rotating the rectangle $\Delta$ by 
$90^{\circ}$, interchanging $a$ and
$b$, because $\cN^{\Delta}(t)=\cN^{\Delta'}(t)$, we also have 
$\cQ(\bbeta(\Delta'),t)=\cQ(\bbeta(\Delta),t)$. Concretely, in this
case we have that counting $(a,b)$-curves on $\PP^1 \times \PP^1$ is the
same as counting $(b,a)$-curves. We thus see that we also have
\[  [t^\delta] \cQ(\bbeta(\Delta),t) = A(\delta) \cdot (2ab) + L(\delta) \cdot (2a+2b) + H(\delta) \cdot a + C(\delta),\]
and since $a \neq b,$ we must have that $H(\delta)=0$.
\end{proof}

Putting together Lemma \ref{lem:h-0} with Corollaries \ref{cor:wH}
and \ref{cor:reversals-Q}, we obtain the following.

\begin{cor}\label{cor:Q-Delta-delta}
Suppose $\Delta$ is an $h$-transverse polygon and all of its edges have 
length at least $\delta.$ Let $v'_i$ be the number of internal vertices of 
determinant $i$. Then
\begin{align}
	Q^{\Delta, \delta} =& A(\delta) \cdot \area(\Delta) + L(\delta) \cdot \LL(\Delta) + D(\delta) \cdot \idet(\Delta) \label{equ:Q-Delta-delta}  \\
&+ C(\delta) + \DQ(\tdet(\Delta), \delta) +\DQ(\bdet(\Delta), \delta) + \sum_{i=1}^{\delta} b_{\delta,i} v'_i. \nonumber 
\end{align}
\end{cor}

Applying Lemma \ref{lem:h-0} also yields a simpler formula for $L(\delta).$

\begin{cor}
  \[ L(\delta) = \frac{1}{2} \sum_{\Gamma: \delta(\Gamma) = \delta} \mu(\Gamma) \eta_0(\Gamma).\]
\end{cor}

We now define a ``correction term'' which will arise in the non-Gorenstein
setting, as well as a more convenient variant of the coefficient $C(\delta)$.

\begin{notn}\label{notn:COR}
For any $\delta \ge 1$,
set
\[ \widetilde{C}(\delta) := C(\delta)-4D(\delta)-4 b_{\delta,1}.\]
Given also $p \geq 1$, set
$$\COR(p,\delta)= 
(2-p) D(\delta) + \DQ(p, \delta) + 2 b_{\delta,1} - b_{\delta, p}
-\frac{1}{6} \widetilde{C}(\delta) \frac{(p-1)(p-2)}{p}.$$
Finally, set $\COR(0,\delta)=0$ for any $\delta \geq 1$.
\end{notn}

We see that the correction term vanishes in the Gorenstein case.

\begin{prop}\label{prop:correction} If $p \leq 2$, then $\COR(p, \delta)=0$.
\end{prop}

\begin{proof} First observe that when $p=1,2$ the last term of 
$\COR(p,\delta)$ clearly vanishes, so it is enough to see that
$$(2-p) D(\delta) + \DQ(p, \delta) + 2 b_{\delta,1} - b_{\delta, p}=0.$$

Now, consider $\Delta$ with $\br_0 = (1, 1, \dots, 1, 0, \dots, 0)$,
$\bl_0 = (0, 0, \dots, 0),$ and $d^t=0$, and having side lengths all at 
least $\delta$. Then applying Corollary \ref{cor:Q-Delta-delta} we find
\begin{align*}
  Q^{\Delta, \delta} =& A(\delta) \cdot \area(\Delta) + L(\delta) \cdot \LL(\Delta) + D(\delta) \cdot 1  + C(\delta) + \DQ(1, \delta) + b_{\delta,1}. \\
\intertext{ Rotating $\Delta,$ as in the proof of Lemma \ref{lem:h-0}, we get}
 Q^{\Delta, \delta} =& A(\delta) \cdot \area(\Delta) + L(\delta) \cdot \LL(\Delta) + D(\delta) \cdot 0  + C(\delta).
\end{align*}
Comparing these two formulas, we get 
$$D(\delta)+\DQ(1, \delta) + b_{\delta,1}=0,$$
yielding the desired statement for $p=1$.

Next consider $\Delta$ with 
$\br_0 = (1, 1, \dots, 1, 0, \dots, 0)$ and 
$\bl_0 = (-1,\cdots, -1, 0, \dots, 0)$, again with $d^t=0$ and having side 
lengths at least
$\delta$. Then Corollary \ref{cor:Q-Delta-delta} gives
\begin{align*}
  Q^{\Delta, \delta} =& A(\delta) \cdot \area(\Delta) + L(\delta) \cdot \LL(\Delta) + D(\delta) \cdot 2  + C(\delta) + \DQ(2, \delta) + b_{\delta,1} \cdot 2. \\
\intertext{ Rotating $\Delta,$ we get}
Q^{\Delta, \delta} =& A(\delta) \cdot \area(\Delta) + L(\delta) \cdot \LL(\Delta) + D(\delta) \cdot 2  + C(\delta) + b_{\delta,2} \cdot 1.
\end{align*}
This yields
$$\DQ(2, \delta) + 2 b_{\delta,1} - b_{\delta,2}=0,$$
and the proposition.
\end{proof}

\begin{notn}\label{notn:det} Given a lattice polygon $\Delta$, let
$\det(\Delta)$ denote the sum of the determinants of the vertices of 
$\Delta$.
\end{notn}

The following theorem is then our main universality result. In view
of Proposition \ref{prop:correction} and Proposition \ref{prop:toric-2}, 
we see that Theorem \ref{thm:main} is just the
special case in which $\Delta$ is strongly $h$-transverse.

\begin{thm}\label{thm:main-2}
Suppose $\Delta$ is an $h$-transverse polygon and all of its edges have 
length at least $\delta.$

Let $(Y(\Delta),\sL)$ be the polarized toric surface corresponding to
$\Delta$, with canonical divisor $\sK$, and second Chern class $c_2$. 
Let $S_i$ be the number of singularities $Y(\Delta)$ of index $i+1$,
and set $\tc_2= \det(\Delta)$ and $S=\sum_{i \geq 1} (i+1) S_i$. 
Then
\begin{align*} 
  Q^{\delta}(Y(\Delta),\sL) =& A(\delta) \cdot \sL^2 
  - L(\delta) \cdot (\sL \cdot \sK) 
  + \frac{1}{12}\widetilde{C}(\delta) \sK^2 \\
  & + \left(\frac{1}{12}\widetilde{C}(\delta)+D(\delta)+ b_{\delta,1}\right) \tc_2  - b_{\delta, 1} S
 + \sum_{i=2}^{\delta} b_{\delta,i} S_{i-1} \\
 & + \COR(\tdet(\Delta), \delta) + \COR(\bdet(\Delta), \delta).
\end{align*}
\end{thm}

\begin{proof}
If $v_i$ is the number of vertices of determinant $i$, we first calculate 
that \eqref{equ:Q-Delta-delta} gives us that
 \begin{align*}
	 Q^{\Delta, \delta} =& A(\delta) \cdot \area(\Delta) 
+ L(\delta) \cdot \LL(\Delta) 
+ D(\delta) \cdot \det(\Delta)  \\
  &+ \widetilde{C}(\delta) + \sum_{i=1}^{\delta} b_{\delta,i} v_i + \COR'(\tdet(\Delta), \delta) + \COR'(\bdet(\Delta), \delta),
\end{align*}
where
\[ \COR'(p, \delta) := 
\begin{cases} (2-p) D(\delta) + \DQ(p, \delta) + 2 b_{\delta,1} 
- b_{\delta, p}, \qquad& p \ge 1, \\ 0, \qquad& p = 0.\end{cases}\]
Indeed, if $\tdet(\Delta)=\bdet(\Delta)=0$, this follows from the identities
$\idet(\Delta) = \det(\Delta)-4,$ $v'_1 = v_1-4,$ and $v'_i = v_i$ for 
$i > 1$.
If $\tdet(\Delta)=0$ and $\bdet(\Delta)=p > 0$, then 
$\idet(\Delta) = \det(\Delta)-2-p,$ and if further $p>1$, then 
$v'_1 = v_1-2,$ and $v'_p = v_p -1,$ and 
$v'_i = v_i$ for $i \neq 1,p$, and the desired identity likewise follows.
On the other hand, if $p=1$, then $v'_1=v_1-3$, and $v'_i=v_i$ for $i>1$,
yielding the same formula.
Finally, the cases that $\tdet(\Delta) > 0$ and $\bdet(\Delta)=0$ or that 
$\tdet(\Delta) > 0$ and $\bdet(\Delta) > 0$ are similar.

The desired statement then follows from the above formula for
$Q^{\Delta,\delta}$ together with Propositions \ref{prop:toric-1} 
and \ref{prop:toric-2}.
\end{proof}

We restate our polynomiality result in the Gorenstein case, and record
that it also implies polynomiality for Severi degrees.

\begin{cor}\label{cor:gorenstein}
The universal linear polynomial
\begin{multline*}
  \widehat{T}_\delta(x, y, z, w; s, s_1, \dots, s_{\delta-1} ) \\
  :=  A(\delta) x  - L(\delta) y + \frac{1}{12}\widetilde{C}(\delta) z  
  + \left(\frac{1}{12} \widetilde{C}(\delta)+D(\delta)+b_{\delta,1}\right) w
  - b_{\delta, 1} s
+ \sum_{i=2}^{\delta} b_{\delta,i} s_{i-1}
\end{multline*}
has the property that for any strongly $h$-transverse polygon $\Delta$ whose edges all have length at least $\delta,$
\[ Q^{\delta}(Y(\Delta), \sL) = \widehat{T}_\delta(\sL^2, \sL \cdot \sK, \sK^2, \tc_2; S, S_1, \dots, S_{\delta-1} ).\]
Hence, let
\[ T_\delta(x, y, z, w; s, s_1, \dots, s_{\delta-1} ) := [t^\delta] \exp \left( \sum_{i \ge 1} \widehat{T}_i(x, y, z, w; s, s_1, \dots, s_{i-1} ) t^i \right).\]
Then we have the polynomial expression
\[ N^{\delta}(Y(\Delta), \sL) = T_\delta(\sL^2, \sL \cdot \sK, \sK^2, \tc_2; S, S_1, \dots, S_{\delta-1} ).\]
\end{cor}

What remains is to compare our results with the previously known Theorem
\ref{thm:goettsche} for arbitrary smooth surfaces, 
showing in particular that our results give new combinatorial formulas for 
the coefficients arising therein, and that the new terms we obtain from
singularities have a simple expression in the context of the 
G\"ottsche-Yau-Zaslow formula. In addition to proving the remaining
results from the introduction, in Corollary \ref{cor:b1-b2} we give formulas 
in terms of our coefficients for the power series $B_1(q),B_2(q)$ of the
G\"ottsche-Yau-Zaslow formula.

We begin by proving the asserted agreement 
of our coefficients with those of Corollary \ref{cor:linearity}.

\begin{proof}[Proof of Proposition \ref{prop:coefs-agree}]
Suppose that $Y(\Delta)$ is smooth, so that $c_2=\tc_2$,
$S_i=0$ for all $i$, and $S=0$. In this case, Tzeng's formula and ours
are both linear in the variables $\sL^2, \sL \cdot \sK, \sK^2, c_2$,
with no constant term, and they necessarily agree whenever $Y(\Delta)$
is smooth. As in the proof of Proposition 2.3 of \cite{go3},
if we consider $(\PP^2,\sO(d))$, corresponding to 
$\Delta$ a right triangle with side length $d$, then we have
$(\sL^2,\sL \cdot \sK, \sK^2, c_2)=(d^2,-3d, 9, 3)$, and if we consider
$(\PP^1 \times \PP^1,\sO(a,b))$, corresponding to $\Delta$ a rectangle
with side lengths $a$ and $b$, we have
$(\sL^2,\sL \cdot \sK, \sK^2, c_2)=(2ab,-2(a+b), 8, 4)$. Setting 
$d=\delta,2\delta,3\delta$ and $(a,b)=(\delta,2\delta)$, we obtain
vectors spanning $\RR^4$,
showing that our linear function must agree with Tzeng's.
\end{proof}

Next, following Qviller ((2.7) of \cite{qv1}), we can re-express the 
the G\"ottsche-Yau-Zaslow formula as follows.

\begin{prop}\label{prop:qviller}
Let
\begin{align*}
  t = DG_2(q) =& \sum_{n \ge 1} n \left(\sum_{d|n} d\right) q^n, \\
  g(t) = (DG_2)^{<-1>}(t) =& \sum_{n \ge 1} \left. \frac{d^{n-1}}{dq^{n-1}} \left( \frac{q}{DG_2(q)} \right)^n \right|_{q=0}  \frac{t^n}{n!}.
\end{align*}

In the variables $x, y, z$ and $w$, the G\"ottsche-Yau-Zaslow formula  
becomes the following:
$$\sum_{\delta \geq 0} T_\delta(x,y,z,w) t^\delta = \frac{(t/g(t))^{\frac{z+w}{12}+\frac{x-y}{2}}B_1\left( g(t) \right))^{z} B_2(g(t))^{y}}{(\Delta(g(t))D^2G_2(g(t))/(g(t))^2)^{\frac{z+w}{24}}}.$$
\end{prop}

We now give a couple of useful formulas involving the $g(t)$ defined in
Proposition \ref{prop:qviller}.

\begin{prop}\label{prop:gt} We have $g(t)=t \cA(t)$, and for any
$i \geq 1$, we have
$$\sum_{\delta \geq i} b_{\delta,i} t^{\delta} =\log\cP(g(t)^i).$$
\end{prop}

Recall that 
$$\cA(t)=\exp\left(-\sum_{\delta \geq 1} 2A(\delta) t^{\delta}\right).$$

\begin{proof} If we take the logarithm of both sides in the
formula of Proposition \ref{prop:qviller} and consider the coefficient
of $x$, we obtain from Proposition \ref{prop:coefs-agree} that 
$$\frac{1}{2} \log(t/g(t)) = \sum A(\delta) t^{\delta},$$
which implies that $g(t) = t \cA(t)$. 
We then use \eqref{eq:log-P} to conclude that
$$\sum_{\delta \geq i} b_{\delta,i} t^{\delta} = \log \cP((t\cA(t))^i)
=\log\cP(g(t)^i),$$
as desired.
\end{proof}

We then obtain the following formulas for $B_1(q)$ and $B_2(q)$.

\begin{cor}\label{cor:b1-b2} With notation as in Theorem \ref{thm:goettsche},
we have 
\begin{align*}
B_1(q) =& \left(\cP(q)\right)^{-1} \cdot \exp\left( 
-\sum_{\delta \ge 1} D(\delta) \left(DG_2(q)\right)^\delta \right), 
\text{ and}\\
B_2(q) =& \exp\left( \sum_{\delta \ge 1} (A(\delta)-L(\delta)) \left(DG_2(q)\right)^\delta \right).
\end{align*}
\end{cor}

\begin{proof}
Again taking the logarithm of both sides in the formula of Proposition 
\ref{prop:qviller} and using Proposition \ref{prop:coefs-agree}, if we 
consider the sum of the coefficients of $x$ and $y$ we find
$$\sum_{\delta \ge 1} (A(\delta) -L(\delta)) t^\delta 
=  \log\left( B_2(g(t)) \right),$$
which, recalling that $g(t)$ is inverse to $DG_2(q)$, yields
the desired formula for $B_2$.
Considering instead the $w$ coefficient minus the $z$ coefficient gives
$$\sum_{\delta \ge 1} (D(\delta)+b_{\delta,1}) t^\delta 
= -\log B_1\left( g(t) \right),$$ 
and we obtain the desired formula using the $i=1$ case of 
Proposition \ref{prop:gt}.
\end{proof}

Finally, we reexpress our formulas in the context of the 
G\"ottsche-Yau-Zaslow formulas, giving an explicit description of the
``correction terms'' arising from the singularities we consider, and
completing the proof of the results stated in the introduction.

\begin{proof}[Proof of Corollary \ref{cor:goettsche-corrected}]
We wish to prove that
\begin{multline*}
   \sum_{\delta \geq 0}T_\delta(x, y, z, w; s, s_1, \dots, s_{\delta-1})(DG_2(q))^\delta \\
  =  \frac{(DG_2(q)/q)^{\frac{z+w}{12}+\frac{x-y}{2}}B_1(q)^{z}B_2(q)^{y}}{(\Delta(q)D^2G_2(q)/q^2)^{\frac{z+w}{24}}} \cP(q)^{-s} \prod_{i \ge 2} \cP\left(q^i \right)^{s_{i-1}}.
\end{multline*}

Now, the formula for 
$\widehat{T}_{\delta}(x,y,z,w;s,s_1,\dots,s_{\delta-1})$ given in
Corollary \ref{cor:gorenstein}, together with Proposition 
\ref{prop:coefs-agree}, gives us that
\begin{align*}
\sum_{\delta \geq 0} T_\delta(x,y,z,w;& s, s_1, \dots, s_{\delta-1}) t^\delta 
\\ 
= & \left(\sum_{\delta \geq 0} T_{\delta}(x,y,z,w) t^{\delta} \right)
\exp\left(-s \sum_{\delta \geq 1} b_{\delta,1} t^{\delta}\right) 
\cdot \prod_{i\geq 2} 
\exp\left(s_{i-1} \sum_{\delta \geq i} b_{\delta,i} t^{\delta}\right) \\
=& \frac{(t/g(t))^{\frac{z+w}{12}+\frac{x-y}{2}}B_1\left( g(t) \right))^{z} B_2(g(t))^{y}}{(\Delta(g(t))D^2G_2(g(t))/(g(t))^2)^{\frac{z+w}{24}}} \cP(g(t))^{-s} \prod_{i \ge 2} \cP\left( \left( g(t) \right)^i \right)^{s_{i-1}},
\end{align*}
with the final expression following from Propositions \ref{prop:qviller}
and \ref{prop:gt}.
The desired statement then follows.
\end{proof}

Note that although we have only proved our enumerative results for Gorenstein 
$h$-transverse toric surfaces, the formula of Corollary 
\ref{cor:goettsche-corrected} 
is unconditional, as it is simply a formal statement on the polynomials 
$T_\delta(x, y, z, w; s, s_1, \dots, s_{\delta-1})$.
Of course, we have only proved that the latter calculate Severi degrees
in the more restrictive setting.
Also, although the product in this formula contains 
infinitely many power series, only finitely many contribute in any given
degree.

\begin{rem}\label{rem:gorenstein} 
We observe that nearly the only case in which the formula of
Proposition \ref{prop:toric-2} for $\sK^2$ yields an integer is the Gorenstein
case. Thus, in the non-Gorenstein case we will likely have to find an
appropriate substitute for $\sK^2$. At the same time, as the singularities
become more complicated, one might expect that a generalized formula has
to take into account more than just the index; in particular, the number
of blowups in a minimal desingularization is no longer determined by the
index. 
\end{rem}

\begin{rem}\label{rem:ampleness}
The Kool-Shende-Thomas proof of universality \cite{k-s-t1}
yields the statement that $N^{\delta}(Y,\sL)$ agrees with
the universal polynomial whenever $\sL$ is $\delta$-very-ample.
We see by Proposition \ref{prop:toric-3} that for the toric surfaces
we consider, $\delta$-very-ampleness implies that all sides of $\Delta$
have length at least $\delta$, so the threshold for polynomiality in our
result is as strong as that of Kool-Shende-Thomas. 

In the special case of $\PP^2$, Kleiman-Shende \cite{k-s1} have proved 
a better (and expected to be 
sharp) threshold of $\lceil \delta/2\rceil+1$, and in fact their
results apply also to Hirzebruch surfaces and classical del Pezzo surfaces,
so at least for some surfaces there is room for further improvement in our 
techniques. In
fact, this sharper threshold is not at all obvious in our combinatorial
setting, and suggests the existence of deeper structure which, if 
understood, could also lead to further simplifications of our formulas.
\end{rem}

\section{Examples of small cogenus}\label{sec:examples}

We conclude by giving some examples, including full calculation of the
coefficients in Theorem \ref{thm:main-2} for $\delta=1,2$ and a look at
a non-Gorenstein example.

\begin{sidewaystable}
\centering
\vspace*{15cm}
\begin{tabular}{c|c|c|c|c|c|c|c|c|c|c|c|c|}
$\Gamma$ &
$\delta(\Gamma)$ & $\ell(\Gamma)$ 
& $\mu(\Gamma)$ 
& $\epsilon_0(\Gamma)$ & $\epsilon_1(\Gamma)$ & $\lambda(\Gamma)$ & $\olam(\Gamma)$ 
& $\Phi(\Gamma, \bbeta)$
& $\zeta^0(\Gamma)$ & $\zeta^1(\Gamma)$ & $\zeta^2(\Gamma)$ & $\eta_0(\Gamma)$
\\
\hline \hline
&&&&&&&&&&&&\\[-.1in]
\begin{picture}(95,8)(-10,-4)\setlength{\unitlength}{2.5pt}\thicklines
\oo
\put(5,2){\makebox(0,0){$\scriptstyle 2$}}
\Eeee
\end{picture}
& 1 & 1 & 4 & 0 & 0 & (2) & (1) &  $\beta_1-1$ 
& 1 & 0 & 0 & $-1$
\\[.15in]
\begin{picture}(95,8)(-10,-4)\setlength{\unitlength}{2.5pt}\thicklines
\ooo
\qbezier(0.8,0.6)(10,5)(19.2,0.6)
\put(10,3.5){\makebox(0,0){$\scriptstyle 1$}}
\end{picture}
& 1 & 2 & 1 & 1 & 1 & (1,1) & (1,1) &  $\beta_1+\beta_2$
& 2 & 1 & 0 & 0
\\
&&&&&&&&&&&&\\[-.1in]
\hline \hline
&&&&&&&&&&&&\\[-.1in]
\begin{picture}(95,8)(-10,-4)\setlength{\unitlength}{2.5pt}\thicklines
\oo
\put(5,2){\makebox(0,0){$\scriptstyle 3$}}
\Eeee
\end{picture}
& 2 & 1 & 9 & 0 & 0 & (3) & (2) &  $\beta_1-2$ 
& 1 & 0 & 0 & -2
\\[.15in]
\begin{picture}(95,8)(-10,-4)\setlength{\unitlength}{2.5pt}\thicklines
\oo
\put(5,3.5){\makebox(0,0){$\scriptstyle 2$}}
\put(5,-3.5){\makebox(0,0){$\scriptstyle 2$}}
\qbezier(0.8,0.6)(5,2)(9.2,0.6)
\qbezier(0.8,-0.6)(5,-2)(9.2,-0.6)
\end{picture}
& 2 & 1 & 16 & 0 & 0 & (4) & (2) & $-\frac{3}{2}\beta_1+\frac{5}{2}$ 
& -$\frac{3}{2}$ & 0 & 0 & $\frac{5}{2}$
\\[.15in]
\begin{picture}(95,8)(-10,-4)\setlength{\unitlength}{2.5pt}\thicklines
\ooo
\qbezier(0.8,0.6)(10,4)(19.2,0.6)
\qbezier(0.8,-0.6)(10,-4)(19.2,-0.6)
\put(10,3.5){\makebox(0,0){$\scriptstyle 1$}}
\put(10,-3.5){\makebox(0,0){$\scriptstyle 1$}}
\end{picture}
& 2 & 2 & 1 & 1 & 1 & (2,2) & (2,2) & $-\frac{3}{2}\beta_1-\frac{3}{2}\beta_2+1$ 
& -3 & -$\frac{3}{2}$ & 0 & 1
\\[.15in]
\begin{picture}(95,8)(-10,-4)\setlength{\unitlength}{2.5pt}\thicklines
\ooo
\qbezier(0.8,0.6)(10,4)(19.2,0.6)
\put(5,-2){\makebox(0,0){$\scriptstyle 2$}}
\Eeee
\put(10,3.5){\makebox(0,0){$\scriptstyle 1$}}
\end{picture}
& 2 & 2 & 4 & 0 & 1 & (3,1) & (2,1) & $-2\beta_1-\beta_2+2$
& -3 &  -1 & 0 & 2 
\\[.15in]
\begin{picture}(95,8)(-10,-4)\setlength{\unitlength}{2.5pt}\thicklines
\ooo
\qbezier(0.8,0.6)(10,4)(19.2,0.6)
\put(15,-2){\makebox(0,0){$\scriptstyle 2$}}
\eEee
\put(10,3.5){\makebox(0,0){$\scriptstyle 1$}}
\end{picture}
& 2 & 2 & 4 & 1 & 0 & (1,3) & (1,2) & $-\beta_1-2\beta_2+2$
& -3 & -2 & 0 & 2 
\\[.15in]
\begin{picture}(95,8)(-10,-4)\setlength{\unitlength}{2.5pt}\thicklines
\oooo
\qbezier(0.8,0.6)(15,6)(29.2,0.6)
\put(15,4.5){\makebox(0,0){$\scriptstyle 1$}}
\end{picture}
& 2 & 3 & 1 & 1 & 1 & (1,1,1) & (1,1,1) & $\beta_1+\beta_2+\beta_3$ 
& 3 & 3 & 1 & 0 
\\[.15in]
\begin{picture}(95,8)(-10,-4)\setlength{\unitlength}{2.5pt}\thicklines
\oooo
\qbezier(0.8,0.6)(10,5)(19.2,0.6)
\qbezier(10.8,0.6)(20,5)(29.2,0.6)
\put(9,4){\makebox(0,0){$\scriptstyle 1$}}
\put(21,4){\makebox(0,0){$\scriptstyle 1$}}
\end{picture}
& 2 & 3 & 1 & 1 & 1 & (1,2,1) & (1,2,1) & $-\beta_1-\beta_2-\beta_3$
& -3 & -3 & -1 & 0 
\\[-.05in]
\end{tabular}
\bigskip
\bigskip
\caption{The templates with $\delta(\Gamma) \le 2$.}
\label{tab:templates}
\end{sidewaystable}

\begin{ex}[Examples of $\delta \le 2$] We compute $A(\delta), L(\delta), C(\delta)$ and $D(\delta)$ for $\delta \le 2$ from the definitions, using
Table \ref{tab:templates}, which was likewise computed from the definitions.
  
 For $\delta = 1:$
\begin{align*}
  A(1) =& \frac{1}{2} \left( 4 \cdot 1 + 1 \cdot 2\right) = 3, \\
  L(1) =& -\frac{1}{2} \left((4 \cdot 1 \cdot (1-0-0) +1 \cdot 2 \cdot (2-1-1) \right) = -2, \\
  D(1) = & -\left( 4 \cdot (0 + 0 \cdot (1-0)) + 1 \cdot (0 + 1 \cdot (1-1)) \right)  = 0, \\
  C(1) = & - \left( 4 \cdot (-1) \cdot (1-0-0) + 1 \cdot 0 \cdot (2-1-1) \right) = 4.
\end{align*}
For $\delta =2:$
\begin{align*}
  A(2) =& \frac{1}{2} \left(9 \cdot 1 + 16 \cdot \left( -\frac{3}{2} \right) + 1 \cdot (-3) + 4 \cdot (-3) + 4 \cdot (-3) + 1 \cdot 3 + 1 \cdot (-3)\right) = -21, \\
  L(2) =& -\frac{1}{2} \begin{pmatrix}
    9 \cdot 1 \cdot (1-0-0)+ 16 \cdot \left( -\frac{3}{2} \right) \cdot (1-0-0) + 1 \cdot (-3) \cdot (2-1-1) \\ + 4 \cdot (-3) \cdot (2-0-1) + 4 \cdot (-3) \cdot (2-1-0) \\
    + 1 \cdot 3 \cdot (3-1-1)+ 1 \cdot (-3)\cdot (3-1-1)
  \end{pmatrix} = \frac{39}{2}, \\
  D(2) =& - \begin{pmatrix}
    9 \cdot (0 + 0 \cdot (1-0)) + 16 \cdot ( 0 + 0 \cdot (1-0)) + 1 \cdot \left( 0 - \frac{3}{2} (1-1) \right)\\ 
    + 4 \cdot (0 - 1 \cdot (1-0)) + 4 \cdot (0 - 2 \cdot (1-1)) \\
    + 1 \cdot (1 + 3 \cdot (1-1)) + 1 \cdot (-1 - 3 \cdot (1-1)) 
  \end{pmatrix}=4, \\
   C(2) =& -\begin{pmatrix}
    9 \cdot (-2) \cdot (1-0-0)+ 16 \cdot \frac{5}{2} \cdot (1-0-0) \\
   + 1 \cdot 1 \cdot (2-1-1) + 4 \cdot 2 \cdot (2-0-1) \\
   + 4 \cdot 2 \cdot (2-1-0) + 1 \cdot 0 \cdot (3-1-1) \\
   + 1 \cdot 0 \cdot (3-1-1)
  \end{pmatrix} = -38.
\end{align*}

Since \[ \cA(t) = \exp\left( -\sum_{\delta \ge 1}  2 A(\delta) t^{\delta} \right),\]
we have
\[ [t^\delta] \cA(t) = \sum_{i \ge 1} \frac{(-1)^i}{i!} \sum_{(\delta_1,\dots, \delta_i)} \prod_{j=1}^i (2 A(\delta_i)),\]
where the summation is over all $i$-compositions of $\delta.$ Hence,
\begin{align*} [t^0] \cA(t) =& 1, \qquad [t^1] \cA(t) = \frac{-1}{1} (2 A(1)) = -6, \\
  [t^2] \cA(t) =& \frac{-1}{1} (2 A(2)) + \frac{1}{2} (2 A(1) \cdot 2A(1)) = 42 +18 = 60.
\end{align*}
So
\[ \cA(t) = 1 - 6 t + 60 t^2 + \cdots , \qquad t \cA(t) = t - 6t^2 + 60t^3 + \cdots . \]
Note that the above formula for $t\cA(t)$ agrees with $g(t)$, as asserted
by Proposition \ref{prop:gt}.
\begin{multline*}
 \sum_{n \ge 1} \frac{\sum_{d | n} d}{n} \cdot (t \cA(t))^{n i} 
= (t - 6t^2 + 60 t^3 + \cdots)^i + \frac{3}{2} (t - 6t^2 + 60t^3 + \cdots)^{2i} \\+ \frac{4}{3} (t - 6t^2 + 60t^3 + \cdots)^{3i} + \cdots .
\end{multline*}
Hence,
$$\begin{array}{lll}  b_{1, 1} = 1,\qquad & b_{2, 1} = - \frac{9}{2},\qquad & b_{3,1} = \frac{130}{3} , \\
  b_{1,2}=0, & b_{2,2} = 1, & b_{3,2} = -12, \\ 
  b_{1,3}= 0, & b_{2,3}=0, & b_{3,3} = 1.
\end{array}$$

We summarize the main coefficients for $\delta=1,2$ in the below table.
\[\begin{tabular}{|c||c|c|c|c|c|c|c|c|}
  \hline
  &&&&&&&& \\[-.14in]
  $\delta$ & $A(\delta)$ & $L(\delta)$ & $D(\delta)$ & $C(\delta)$ & $b_{\delta,1}$ & $b_{\delta,2}$ & $\widetilde{C}(\delta)$ & $\frac{1}{12}\widetilde{C}(\delta)+D(\delta)+b_{\delta,1}$ \\[.02in]
  \hline
  &&&&&&&& \\[-.15in]
  \hline
  $1$ & $3$ & $-2$ & $0$ & $4$ & $1$ & $0$ & $0$ & $1$ \\
  \hline
  &&&&&&&& \\[-.14in]
  $2$ & $-21$ & $\frac{39}{2}$ & $4$ & $-38$ & $-\frac{9}{2}$ & $1$ & $-36$ & $-\frac{7}{2}$ \\[.02in]
  \hline
\end{tabular}\]
\end{ex}

\begin{ex} From the previous example, we can compute our universal linear
functions in the cases $\delta=1,2$.
Suppose $\Delta$ is an $h$-transverse polygon and all of its edges have length at least $\delta.$ Let $v'_i$ be the number of internal vertices of determinant $i.$ Then Corollary \ref{cor:Q-Delta-delta} yields
\begin{align*}
  Q^{\Delta, 1} =& 3 \cdot \area(\Delta) - 2 \cdot \LL(\Delta) + 4  + \DQ(\tdet(\Delta), 1) +\DQ(\bdet(\Delta), 1)  + v'_1, \\
  Q^{\Delta, 2} =& -21 \cdot \area(\Delta) + \frac{39}{2} \cdot \LL(\Delta) + 4 \cdot \idet(\Delta) - 38 \\
  & + \DQ(\tdet(\Delta), 2) +\DQ(\bdet(\Delta), 2) - \frac{9}{2} v'_1 + v'_2.
\end{align*}
 
If further $\Delta$ is strongly $h$-transverse, then
with notation as in Theorem \ref{thm:main-2}, we have
\begin{align*}
  Q^{1}(Y(\Delta),\sL) =& 3 \cdot \sL^2 + 2 \cdot (\sL \cdot \sK) + 0 \cdot \sK^2  + 1 \cdot \tc_2 - 1\cdot S, \\
  Q^{2}(Y(\Delta),\sL) =& -21 \cdot \sL^2 -\frac{39}{2} \cdot (\sL \cdot \sK) - 3 \cdot \sK^2 - \frac{7}{2} \cdot \tc_2 + \frac{9}{2} S + 1 \cdot S_1.
\end{align*}

If we apply this to the smooth cases, in which $S=S_i=0$ and $\tc_2=c_2,$ we get
\begin{align*}
  Q^{1}(Y(\Delta),\sL) =& 3 \cdot \sL^2 + 2 \cdot (\sL \cdot \sK) + 0 \cdot \sK^2  + 1 \cdot c_2, \\
  Q^{2}(Y(\Delta),\sL) =& -21 \cdot \sL^2 -\frac{39}{2} \cdot (\sL \cdot \sK) 
 - 3 \cdot \sK^2 - \frac{7}{2} \cdot c_2.
\end{align*}
These agree with previously obtained formulas; see, e.g., the formulas at
the beginning of Kleiman-Piene \cite{k-p2} (in their notation, $a_1$ and
$a_2/2$).
\end{ex}

\begin{ex} We conclude with an example of what happens in the non-Gorenstein
case, calculating the resulting correction terms explicitly.
  Consider $\Delta$ with $\tdet(\Delta)=3$ and $\bdet(\Delta)=0.$ Applying Lemma \ref{lem:DQ-pgedelta}, we get for $p \ge 2,$
\begin{align*}
\DQ(p,1) & = -  \sum_{\Gamma: \epsilon_0(\Gamma)=1, \delta(\Gamma)=1} \mu(\Gamma) \left( p \zeta^1(\Gamma) + \eta_0(\Gamma) \right) = -p \\
\DQ(p,2) & = -   \sum_{\Gamma: \epsilon_0(\Gamma)=1, \delta(\Gamma)=2} \mu(\Gamma)\left(p \zeta^1(\Gamma) + \eta_0(\Gamma) \right) = \frac{19}{2}p -9. 
\end{align*}
Then Corollary \ref{cor:Q-Delta-delta} gives us
\begin{align*}
  Q^{\Delta, 1} =& 3 \cdot \area(\Delta) - 2 \cdot \LL(\Delta) + 4 -3  + v'_1 \\
  Q^{\Delta, 2} =& -21 \cdot \area(\Delta) + \frac{39}{2} \cdot \LL(\Delta) + 4 \cdot \idet(\Delta)  -38 + \frac{39}{2} - \frac{9}{2} v'_1 + v'_2. 
\end{align*}
If we express the formula using Theorem \ref{thm:main-2}, we instead find
\begin{align*}
  Q^{\Delta, 1} =& 3 \cdot \sL^2 + 2 \cdot (\sL \cdot \sK) + 0 \cdot \sK^2 
+1 \cdot \tc_2  - 1 \cdot S + \COR(3,1) \\
 = & 3 \cdot \sL^2 + 2 \cdot (\sL \cdot \sK) + 0 \cdot \sK^2 
+1 \cdot \tc_2  - 1 \cdot S - 1 \\
  Q^{\Delta,2} =& -21 \cdot \sL^2 -\frac{39}{2} \cdot (\sL \cdot \sK) - 3 \cdot \sK^2 - \frac{7}{2} \cdot \tc_2 + \frac{9}{2} S + 1 \cdot S_1+\COR(3,2)\\
   =& -21 \cdot \sL^2 -\frac{39}{2} \cdot (\sL \cdot \sK) - 3 \cdot \sK^2 - \frac{7}{2} \cdot \tc_2 + \frac{9}{2} S + 1 \cdot S_1+\frac{21}{2}.
\end{align*}
\end{ex}

\appendix
\section{Calculations on toric surfaces}\label{app:toric}

In this appendix, we collect some explicit calculations on toric surfaces
which are presumably well known, but for which we do not have references.
Throughout this appendix, $\Delta$ denotes an $h$-transverse polygon,
and $(Y(\Delta),\sL)$ the corresponding polarized toric variety with
canonical divisor $\sK$.
We denote by $v_i$ the number of vertices of $\Delta$ having determinant
$i$, and $S_i$ the number of singularities of $Y(\Delta)$ of index
$i+1$.

\begin{prop}\label{prop:toric-1} 
For any $v \in \Delta$ with $\det(v)>1$, we have that $\det(v)$ is 
the index of the singularity corresponding to $v$, so 
$$v_i=S_{i-1} \text{ for } i>1.$$

In addition, $Y(\Delta)$ is Gorenstein if and only if it has only rational
double points, if and only if
$\tdet(\Delta), \bdet(\Delta) \in \{0,1,2\}$, and in this case every
$v$ with $\det(v)>1$ corresponds to a singularity whose minimal resolution
introduces a chain of $\det(v)-1$ rational curves.
\end{prop}

\begin{proof}
First, for any vertex $v \in \Delta$, we have that $v$ corresponds to 
a singularity of $Y(\Delta)$ exactly when $\det(v)>1$, and that the
corresponding singularity is a cyclic quotient singularity of index 
$\det(v)$; see \S 10.1 of \cite{c-l-s1}.

Now, if $v$ is an internal
vertex, the adjacent outward normals are either of the form 
$(1,a),(1,b)$ with $a<b$ positive integers, or $(-1,a),(-1,b)$ with $a>b$ 
positive integers. According to \S 10.1 of \cite{c-l-s1}, we should perform
a $GL_2(\ZZ)$ operation on these normals to put them in the form
$(0,1),(d,-k)$ where $d>0$, $0 \leq k < d$, and $\gcd(k,d)=1$.
In the first case, we find that $d=b-a$ and $k=d-1$, while in
the second case, we find that $d=a-b$ and $k=d-1$. Thus, according to
Proposition 10.1.6 of \cite{c-l-s1}, the corresponding singularity is
Gorenstein, and by \S 10.4,
a minimal resolution of such singularities has a chain of $d-1$ rational
curves lying over the singularity, as claimed.

Proposition 10.4.11 (together with Proposition 10.1.6) of \cite{c-l-s1} says 
that being Gorenstein is equivalent to having only rational double points
for singularities. For the next assertion, it remains to consider 
non-internal vertices. If $\tdet(\Delta) = 0$ (i.e., $d^t \neq 0$), then 
there are two vertices at the top of $\Delta$, which each have determinant 
$1$, and hence do not yield any singularities. If $\tdet(\Delta)=d>0$, then 
the adjacent normals at the top are $(1,a),(-1,b)$, and calculating as above 
we find that $d=a+b$ and $k=1$. We see that unless $d=1$, in which case 
there is no corresponding singularity, the only Gorenstein case is when 
$k=1=d-1$, which is to say $d=2$. The calculation is similar for $\bdet$, 
yielding the desired criterion for when $Y(\Delta)$ is Gorenstein.
\end{proof}

\begin{prop}\label{prop:toric-2} 
We have
$$\area(\Delta) = \sL^2, \quad \LL(\Delta)=-\sL \cdot \sK,
\quad \quad \sum_{i \geq 1} v_i = c_2.$$

Further, we can calculate
$$\det(\Delta)=12-\sK^2 + \COR''(\tdet(\Delta)) + \COR''(\bdet(\Delta)),$$
where
\[ \COR''(p) := \begin{cases} \frac{2(p-1)(p-2)}{p}, \qquad& \forall p \ge 1,
\\ 0, \qquad& p = 0.\end{cases}\]
and in the Gorenstein case, this reduces to
$$\det(\Delta)=12-\sK^2.$$
\end{prop}

\begin{proof}
The first two statements are almost the same as Proposition
10.5.6 of \cite{c-l-s1},
with the only difference being that they consider the pullback $\sL'$ of $\sL$
to a desingularization $Y'$ of $Y(\Delta)$. That is, they show that
$\area(\Delta)=(\sL')^2$, and $\LL(\Delta)=-\sL' \cdot \sK_{Y'}$. But
$(\sL')^2=(\sL)^2$ because intersection commutes with pullback, and
$\sL' \cdot \sK_{Y'}=\sL \cdot \sK$ because $\sK_{Y'}$ differs from the
pullback of $\sK$ by a sum of exceptional divisors,
which intersect $\sL'$ trivially because $\sL$ is ample.
We thus conclude the first two formulas. 

The formula for $c_2$ is a special case of the main theorem of
\cite{b-b-f1}. 

Next, we compute $\sK^2$. First, if $\Delta$ is any lattice polygon, using
\S 8.2 and Proposition 6.3.8 \cite{c-l-s1} we conclude that 
$$\sK^2=\sum_{i=0}^{n-1} \left(\frac{1}{d_{i-1,i}}+\frac{1}{d_{i,i+1}}-
\frac{d_i}{d_{i-1,i}d_{i,i+1}}\right),$$
where we number the rays of the normal fan $\rho_0,\dots,\rho_{n-1}$ in
the counterclockwise direction, and define $d_{i,i+1}$ to be the
determinant of the matrix formed by (primitive integer representatives of) 
$\rho_i$ and $\rho_{i+1}$, and $d_i$ to be the determinant of the matrix
formed by $\rho_{i-1}$ and $\rho_{i+1}$. Here indices are taken modulo $n$
as needed. In the case that $\Delta$ is $h$-transverse and $d^t,d^b >0$,
suppose that $\rho_{i_1}=(0,1)$ and $\rho_{i_2}=(0,-1)$. We then see that 
all terms of the above sum vanish except for three each at the
top and bottom of $\Delta$ which involve $\rho_{i_1}$ or $\rho_{i_2}$. 
Let $l_1,\dots,l_M$ and $r_1,\dots,r_M$ denote the (unreordered) left
and right directions of $\Delta$.
Then, if the three terms coming from the top of $\Delta$ are disjoint
from the three coming from the bottom, the top three contribute 
$4-(r_1-l_1)$ 
and the bottom three contribute $4-(l_M-r_M)$,
so the total is 
$$8-(r_1-r_M)+(l_1-l_M)=8-\idet(\Delta)=12-\det(\Delta),$$ 
as desired.
One checks that when $\Delta$ has no internal vertices on the left and/or
right, so that the six aforementioned terms are not disjoint, the same
formula still holds.

Now, if $d^t=0$, the three nonzero terms at the top are replaced by
two nonzero terms, which together contribute $\frac{4}{r_1-l_1}$. This
replaces the contribution of $4-(r_1-l_1)$, and when $d^t=0$ the
formula for the determinant changes by $(r_1-l_1)-2$, and we see that
$$\frac{4}{r_1-l_1}-(4-(r_1-l_1))+((r_1-l_1)-2)=\COR''(\tdet(\Delta)),$$
yielding the desired formula. The calculation at the bottom is the same,
so we obtain the desired formula for $\sK^2$.
\end{proof}

\begin{prop}\label{prop:toric-3} 
If $\sL$ is $d$-very-ample for some $d$, then every edge of 
$\Delta$ has length at least $d$.
\end{prop}

\begin{proof} 
First recall that $d$-very-ample is defined by
the condition that if $Z \subseteq Y$ is any finite subscheme of length
$d+1$, then the restriction map $H^0(Y,\sL) \to H^0(Z,\sL|_Z)$ is surjective.
Now, if $\sL$ is $d$-very-ample, then by definition
it is $d$-very-ample after restricting to any torus orbit closure. Each edge 
of $\Delta$ corresponds to a torus orbit closure isomorphic to $\PP^1$, and
using the explicit descriptions of $\sL$ and of orbit closures given in 
\S 4.2
and Proposition 3.2.7 of \cite{c-l-s1},
one computes that if an edge has length $\ell$, the restriction of $\sL$ to
the corresponding $\PP^1$ has degree $\ell$.
Now, a line bundle of $\PP^1$ is $d$-very-ample if and only it has degree
at least $d$, giving the desired assertion.
\end{proof}

\section*{Index of notation}

\begin{tabular}{llcll}
$N^{\delta}(Y,\sL)$ \quad & Notation \ref{notn:Ndelta} 
& \quad \quad \quad \quad\quad \quad \quad&
$\mu(G)$ & Definition \ref{defn:mudelta} \\
$N^{\Delta,\delta}$ & Notation \ref{notn:NDeltadelta} & \quad &
$\epsilon_i(G)$ & Notation \ref{notn:epsilon} \\
$N_{\bbeta}^{\delta}$ & \eqref{equ:NbetaDelta} & \quad &
$\eta_j(\Gamma)$ & Notation \ref{notn:zetaeta} \\
$\cN(\bbeta,t)$ & \eqref{equ:Nbetat} & \quad &
$\zeta^i(\Gamma)$ & Notation \ref{notn:zetaeta} \\
$\cN^{\Delta}(t)$ & \eqref{equ:NDeltat} & \quad &
$\lambda_j(G)$ & Definition \ref{defn:lambda} \\
$Q^{\delta}(Y,\sL)$ & Notation \ref{notn:Qdelta} & \quad &
$\olam_j(G)$ & Definition \ref{defn:olam} \\
$Q^{\Delta,\delta}$ & Notation \ref{notn:NDeltadelta} & \quad &
$\rho(e)$ & Definition \ref{defn:long-edge} \\
$Q_{\bbeta}^{\delta}$ & \eqref{equ:Qbetat} & \quad &
$\Phi(G,\bbeta)$ & Theorem \ref{thm:linear0} \\
$Q_{\Gamma}(\bbeta)$ & Notation \ref{notn:QGammabbeta} & \quad &
$\Phi_{\bbeta}(G)$ & \eqref{equ:phibeta} \\
$Q_{\delta}(\bbeta)$ & Notation \ref{notn:QGammabbeta} & \quad &
$\Phi_{\bbeta}^s(G)$ & \eqref{equ:phibetas} \\
$\cQ(\bbeta,t)$ & \eqref{equ:Qbetat} & \quad &
$\nBar{\Gamma}$ & Definition \ref{defn:conjugate} \\
$\cQ^{\Delta}(t)$ & \eqref{equ:QDeltat} & \quad &
$\COR(p,\delta)$ & Notation \ref{notn:COR} \\
$\cA(t)$ & Notation \ref{notn:AP} & \quad &
$\DQ(p,\delta)$\quad & Notation \ref{notn:DQ} \\
$\cP(x)$ & Notation \ref{notn:AP} & \quad &
$\bdet(\Delta)$ & Notation \ref{notn:tdet} \\
$[t^{\delta}]$ & Notation \ref{notn:Qdelta} & \quad &
$\tdet(\Delta)$ & Notation \ref{notn:tdet} \\
$b_{\delta,i}$ & \eqref{equ:bdeltai} & \quad &
$\area(\Delta)$ & Corollary \ref{cor:wH} \\
$d^b$ & Notation \ref{notn:dt} & \quad &
$\height(\Delta)$ & Corollary \ref{cor:wH} \\
$d^t$ & Notation \ref{notn:dt} & \quad &
$\idet(\Delta)$ & Corollary \ref{cor:wH} \\
$B_i(q)$ & Theorem \ref{thm:goettsche} & \quad &
$\LL(\Delta)$ & Corollary \ref{cor:wH} \\
$A(\delta)$ & Notation \ref{notn:coefs} & \quad &
$\det(\Delta)$ & Notation \ref{notn:det} \\
$C(\delta)$ & Notation \ref{notn:coefs} & \quad &
$\area(\bbeta)$ & Notation \ref{notn:areaetc} \\
$\widetilde{C}(\delta)$ & Notation \ref{notn:COR} & \quad &
$\height(\bbeta)$ & Notation \ref{notn:areaetc} \\
$D(\delta)$ & Notation \ref{notn:coefs} & \quad &
$\idet(\bbeta)$ & Notation \ref{notn:areaetc} \\
$H(\delta)$ & Notation \ref{notn:coefs} & \quad &
$\LL(\bbeta)$ & Notation \ref{notn:areaetc} \\
$L(\delta)$ & Notation \ref{notn:coefs} & \quad &
$\det(v)$ & Definition \ref{defn:det} \\
$G_{(k)}$ & Definition \ref{defn:shift} & \quad &
$\ext_{\bbeta}(G)$ & Definition \ref{defn:P} \\
$P_{\bbeta}(G)$ & Definition \ref{defn:P} & \quad &
$\maxv(G)$ & Definition \ref{defn:shift} \\
$P_{\bbeta}^s(G)$ & Definition \ref{defn:P} & \quad &
$\minv(G)$ & Definition \ref{defn:shift} \\
$\bbeta(\bd)$ & Notation \ref{notn:divergence} & \quad &
$\Rev(\bs)$ & Definition \ref{defn:Rev} \\
$\bbeta(\Delta)$ & Notation \ref{notn:bbetaDelta} & \quad &
$\ell(G)$ & Definition \ref{defn:shift} \\
$\delta(\bl,\br)$ & Definition \ref{defn:Rev} & \quad & 
$\ell(\Delta)$ & Notation \ref{notn:ellDelta} \\
$\delta(G)$ & Definition \ref{defn:mudelta} & \quad & & \\
\end{tabular}

\bibliographystyle{amsalpha}
\bibliography{gen}

\end{document}

%% file: headers.tex
\usepackage{amssymb,euscript,amsmath, mathrsfs}
\usepackage{bm, rotating}
\makeindex

\newcounter{ENUM}
\newcommand{\itm}{\item}
\newenvironment{ilist}[1][0]{\renewcommand{\theENUM}{\roman{ENUM}}\renewcommand{\itm}{\addtocounter{ENUM}{1}\item[(\theENUM)]}\begin{itemize}\setcounter{ENUM}{#1}}{\end{itemize}}

\newenvironment{alist}[1][0]{\renewcommand{\theENUM}{\alph{ENUM}}\renewcommand{\itm}{\addtocounter{ENUM}{1}\item[(\theENUM)]}\begin{itemize}\setcounter{ENUM}{#1}}{\end{itemize}}

\newcommand{\oo}{\multiput(0,0)(10,0){2}{\circle*{2}}}
\newcommand{\ooo}{\multiput(0,0)(10,0){3}{\circle*{2}}}
\newcommand{\oooo}{\multiput(0,0)(10,0){4}{\circle*{2}}}
\newcommand{\Eeee}{\put(1,0){\line(1,0){8}}}
\newcommand{\eEee}{\put(11,0){\line(1,0){8}}}

\input xy
\xyoption{all}
\CompileMatrices

\def\bd{\bm{d}}
\def\bl{\bm{l}}
\def\br{\bm{r}}
\def\bs{\bm{s}}
\def\bv{\bm{v}}
\def\tc{\tilde{c}}
\def\bbeta{{\bm{\beta}}}

\def\maxv{\mathrm{maxv}}
\def\minv{\mathrm{minv}}
\def\Rev{\mathrm{Rev}}
\def\ext{\mathrm{ext}}
\def\area{\mathrm{area}}
\def\LL{\mathrm{LL}}
\def\height{\mathrm{height}}
\def\idet{\mathrm{idet}}
\def\tdet{\mathrm{tdet}}
\def\bdet{\mathrm{bdet}}
\def\det{\mathrm{det}}
\newcommand{\nBar}[1]{\overline{#1}}
\def\DQ{\mathrm{DiffQ}}

\def\COR{\mathrm{COR}}

\def\ZZ{{\mathbb Z}}

\def\AA{{\mathbb A}}

\def\PP{{\mathbb P}}
\def\QQ{{\mathbb Q}}
\def\RR{{\mathbb R}}

\def\cA{{\mathcal A}}

\def\cN{{\mathcal N}}
\def\cP{{\mathcal P}}
\def\cQ{{\mathcal Q}}

\def\sK{{\mathscr K}}
\def\sL{{\mathscr L}}

\def\sO{{\mathscr O}}

\def\olam{\overline{\lambda}}

\newtheorem{thm}{Theorem}[section]
\newtheorem{prop}[thm]{Proposition}
\newtheorem{lem}[thm]{Lemma}
\newtheorem{cor}[thm]{Corollary}

\theoremstyle{definition}
\newtheorem{defn}[thm]{Definition}

\newtheorem{ex}[thm]{Example}

\theoremstyle{remark}
\newtheorem{notn}[thm]{Notation}
\newtheorem{rem}[thm]{Remark}

\numberwithin{equation}{section}